\documentclass[proof]{sn-jnl}
\usepackage[utf8]{inputenc}
\usepackage[T1]{fontenc}
\usepackage{amsmath}
\usepackage{amsfonts}
\usepackage{amssymb,cleveref,amsthm,xcolor}

\usepackage{xpatch}
\makeatletter
\newcounter{proofpart}
\xpretocmd{\proof}{\setcounter{proofpart}{0}}{}{}
\newcommand{\fact}[1]{%
	\par
	\addvspace{\medskipamount}%
	\stepcounter{proofpart}%
	\noindent\emph{Fact \theproofpart: #1}\par\nobreak\smallskip
	\@afterheading
}
\makeatother

\newcounter{caso}
\xpretocmd{\proof}{\setcounter{caso}{0}}{}{}
\newcommand{\caso}[1]{%
	\par
	\addvspace{\medskipamount}%
	\stepcounter{caso}%
	\noindent\emph{Case \thecaso: #1}\par\nobreak\smallskip
}
\makeatother

\theoremstyle{plain}
\newtheorem{lemma}{Lemma}
\newtheorem{thm}{Theorem}
\newtheorem{prop}{Proposition}
\newtheorem*{theorem*}{Theorem}

\theoremstyle{definition}
\newtheorem{definition}{Definition}

\theoremstyle{remark}
\newtheorem{rmrk}{Remark}
\newtheorem{exmpl}{Example}

\Crefname{thm}{Theorem}{Theorems}
\Crefname{prop}{Proposition}{Propositions}
\Crefname{cor}{Corollary}{Corollaries}

\DeclareMathOperator{\tr}{tr}
\DeclareMathOperator{\imm}{\mathrm{Im}}

\newcommand{\cV}{\mathcal{V}}
\newcommand{\cA}{\mathcal{A}}

\title[Second Variation \& asymptotics]{Operators arising as Second Variation of optimal control problems and their spectral asymptotics}
\author{Stefano Baranzini}
\affil{SISSA, Scuola Internazionale Superiore di Studi Avanzati, Via Bonomea, 265 - 34136 Trieste, Italy} 
\email{sbaranzi@sissa.it}	

\abstract{We compute the asymptotic for the eigenvalues of a particular class of compact operators deeply linked with the second variation of optimal control problems. We characterize this family in terms of a set of finite dimensional data and we apply this results to a particular class of singular extremal to get a nice description of the spectrum of the second variation.}
\keywords{second variation, optimal control, Weyl law, compact operator}

\begin{document}
	\maketitle	
		
	\section*{Introduction}
	
	The main focus of this paper is the study of a particular class of compact operators $K$ on the Hilbert space $L^2([0,1],\mathbb{R}^k)$ with the  standard Hilbert structure. They are characterized by the following properties:
	\begin{itemize}
		\item there exists a finite dimensional subspace of $L^2([0,1],\mathbb{R}^k)$, which we call $\cV$, on which $K$ becomes a self-adjoint operator, i.e. :
		\begin{equation}
			\label{eq: K restricted to V is self adjoint}
			\langle u, K v \rangle  = \langle  K u, v\rangle \quad  \forall\,  u, v \in \cV,
		\end{equation} 
	\item $K$ is an Hilbert-Schmidt operator with an integral kernel of a particular form, namely:
	\begin{equation}
		\label{eq: being volterra type}
		K(v)(t) = \int_0^tV(t,\tau)v(\tau) d\tau, \quad v \in L^2([0,1],\mathbb{R}^k).
	\end{equation}
	Where $V(t,\tau)$ is a matrix whose entries are $L^2$ functions. We call the class of operator satisfying this last condition \emph{Volterra-type} operators.
	\end{itemize}
	
	The main results of this paper are a fairly general study of the asymptotic distribution of the eigenvalues of $K$ when restricted to any subspace $\cV$ which satisfies \cref{eq: K restricted to V is self adjoint} (\Cref{thm: eigenvalues of K}) and a characterization result for operators satisfying the two properties stated above (\Cref{thm: characterization of K}).
	
	The first result is proved in \Cref{section: proof thm 1}. We first restrict ourself to operators $\tilde{K}$ of the form:
		\begin{equation}
		\label{eq: compact part second variation}
		\tilde K(v)(t) = -\int_0^t\sigma (Z_\tau v_\tau,Z_t \cdot) d\tau.
	\end{equation} 
	Here $Z_t$ is an analytic in $t$, $2n\times k $ matrix  and $\sigma$ the standard symplectic form on $\mathbb{R}^{2n}$ (see \Cref{rmrk: symplectic}).  
	A similar asymptotic formula  was proved in \cite[Theorem 1]{determinant}, it was shown that if we consider $\{\lambda_n(\tilde{K})\}_{n \in \mathbb{Z}}$ the decreasing (resp. increasing) arrangement of positive (resp. negative) eigenvalues of $\tilde K$ we have either:
	\begin{equation}
		\label{eq: asymptotyc capacity}
		\lambda_n(\tilde{K}) = \frac{\xi}{\pi n} + O(n^{-5/3}) \quad \text{ or } \quad  \lambda_n(\tilde K) = O(n^{-2}),
	\end{equation}
	for $n \in \mathbb{Z}$ sufficiently large and for some $\xi>0$. The number $\xi$ is called \emph{capacity} and depends only on the matrix $Z_t$ in the definition of $\tilde{K}$. 	
		
	If $\xi=0$, we go further with the expansion in \cref{eq: asymptotyc capacity}. We single out the term giving the principal contribution to the asymptotic representing the quadratic form associated to $\tilde K$ as:
	\begin{equation*}
		Q(v) = \langle v,\tilde K v\rangle = -\int_0^1 \int_0^t \sigma(Z_\tau  v_\tau ,Z_tu_t ) d \tau dt = \sum_{i=1}^{k-1} Q_i(v)+R_k(v).
	\end{equation*}
	The result mentioned above corresponds to  the case $Q_1 \ne 0$, in \Cref{thm: eigenvalues of K} we give the asymptotic for the general case.
	
	From the point of view of geometric control theory \Cref{thm: eigenvalues of K} can be seen as an asymptotic analysis of the spectrum of the second variation for particular classes of singular extremals and a quantitative version of some necessary optimality conditions.
	
	Precise definitions will be given in \Cref{section: control}, standard references on the second variation are \cite[Chapter 20]{bookcontrol} and \cite{ASZ}. For now it is enough to know that the second variation $Q$ of an optimal control problem on a manifold $M$ is a linear operator on $L^2([0,1],\mathbb{R}^k)$ of the following form:
	\begin{equation}
		\label{eq: second variation}
		\langle Q v, u \rangle = -\int_0^1 \langle H_t v_t,u_t \rangle -\int_0^1 \int_0^t \sigma(Z_\tau  v_\tau ,Z_tu_t ) d \tau dt,
	\end{equation}
	where $H_t$ is a symmetric $k\times k$  matrix, $\sigma$ is the standard symplectic form on $T_\eta T^*M$ and $Z_t: \mathbb{R}^k \to T_\eta(T^*M)$ is a linear map with values in the tangent space to a fixed point $\eta \in T^*M$.

	For totally singular extremal, the matrix $H_t$ appearing in \cref{eq: second variation} is identically zero and the second variation reduces to an operator of the same form as in \cref{eq: compact part second variation}.

	In \Cref{section: proof thm 2} we prove \Cref{thm: characterization of K}. We first show that any $K$ satisfying \cref{eq: K restricted to V is self adjoint,eq: being volterra type} it is completely determined by its (\emph{finite rank}) skew-symmetric part $\cA$ and can always be represented as in \cref{eq: compact part second variation}. Then we relate the \emph{capacity} of $K$ to the spectrum of $\cA$.  
	

	In \Cref{section: control} we recall some  basic notions from control theory and we reformulate \Cref{thm: characterization of K} in a more control theoretic fashion, and use it to characterize the operators coming form the \emph{second variation} of an optimal control problem. Moreover we give a geometric interpretation of the capacity $\xi$ appearing in \cref{eq: asymptotyc capacity} in terms of the Hessian of the maximized Hamiltonian coming from Pontryagin Maximum Principle. 		 
		
 \section{Overview of the main results}
		We begin this section recalling some general facts about the spectrum of compact operators, then we fix some notation and give a precise statement of the main results. Given a compact self-adjoint operator $K$ on an Hilbert space $\mathcal{H}$, we can define a quadratic form  setting $Q (v) = \langle v,K(v)\rangle$. The eigenvalues of $Q$ are by definition those of $K$ and we will denote $\Sigma_{\pm}(Q)$ the positive and negative parts of the spectrum of $Q$. 

By the standard spectral theory of compact operators (see \cite{functionalanalysis}) the non zero eigenvalues of $K$ are either finite or accumulate at zero and their multiplicity is finite.  Consider the positive part of the spectrum of $Q$, $\Sigma_{+}(Q)$ and $\lambda \in \Sigma_+(Q)$. Denote by $m_\lambda$ the multiplicity of the eigenvalue $\lambda$. We can introduce a monotone non increasing sequence $\{\lambda_n\}_{n \in \mathbb{N}}$ indexing the eigenvalues of $K$, requiring that the cardinality of the set $\{\lambda_n = \lambda\} = m_\lambda$ for every $\lambda \in \Sigma_+(Q)$.

This will be called the monotone arrangement of $\Sigma_+(Q)$. We can  perform the same construction indexing by $-n$, $n \in \mathbb{N}$, the negative part of the spectrum $\Sigma_-(Q)$. This time we require that the sequence $\{\lambda_{-n}\}_{n \in \mathbb{N}}$ is non decreasing. Provided that $\Sigma_\pm(Q)$ are both infinite, we obtain a sequence $\{\lambda_n\}_{n \in \mathbb{Z}}$.
\begin{definition}
	\label{def: capacity}
	Let $Q$ be a quadratic form $Q$ on a Hilbert space $\mathcal{H}$ and $j \in \mathbb{N}$
	\begin{itemize}
		\item if $j$ is odd, $Q$ has $j-$capacity $\xi>0$ with reminder of order $\nu>0$ if $\Sigma_+(Q)$ and $\Sigma_{-}(Q)$ are both infinite and:
	\begin{equation*}
		\lambda_n = \frac{\xi}{(\pi n)^j} + O(n^{-\nu-j}) \quad \text{ as }\quad n \to \pm \infty,
	\end{equation*}
	\item if $j$ is even, $Q$ has $j-$capacity $(\xi_+,\xi_-)$ of order $\nu >0$ if both  $\Sigma_+(Q)$ and $\Sigma_{-}(Q)$ are infinite and:
	\begin{equation*}
		\begin{split}
		\lambda_n = \frac{\xi_+}{(\pi n)^j} + O(n^{-\nu-j}) \quad \text{ as }\quad n \to + \infty,\\
		\lambda_n = \frac{\xi_-}{(\pi n)^j} + O(n^{-\nu-j}) \quad \text{ as }\quad n \to - \infty,	
		\end{split}
	\end{equation*}
	where $\xi_\pm\ge0$ or if at least one between $\Sigma_+(Q)$ and $\Sigma_-(Q)$ is infinite and the relative monotone arrangement satisfies the corresponding asymptotic relation;
	\item if the spectrum is finite or $\lambda_n = O(n^{-\nu})$ as $n \to \pm \infty$ for any $\nu >0$, we say that $Q$ has $\infty-$capacity.
	\end{itemize} 
\end{definition}

The behaviour of the sequence $\{\lambda_n\}_{n \in \mathbb{Z}}$ is closely related to the following counting functions:
\begin{equation*}
	C^+_j(n) =\# \{l \in \mathbb{N} : 0< \frac{1}{\sqrt[j]{\lambda_l}}<n\} \quad 	C^-_j(n) =\# \{l \in \mathbb{N} : -n> \frac{-1}{\sqrt[j]{\vert\lambda_{-l}\vert}}>0\}
\end{equation*}

The requirement of \cref{def: capacity} for the $j-$capacity can be translated into the following asymptotic for the functions $C^{\pm}_j(n)$:
\begin{equation*}
	C^\pm_j(n) = \frac{\xi_\pm}{\pi} n + O(n^{1-\nu}) \quad \text{ as } \quad n \to \pm \infty
\end{equation*}

We illustrate here some of the properties of the $j-$capacity. The proofs are given in \cref{section: proof thm 1}, \Cref{prop capacity}. Without loss of generality we state the properties for the positive part of the spectrum, analogue results hold for the negative one. 
\begin{itemize}
	\item \textit{(Homogeneity)} if $Q_1$ and $Q_2$ are quadratic forms on two Hilbert spaces $\mathcal{H}_1$ and $\mathcal{H}_2$ of $j-$capacity $\xi_1$ and $\xi_2$ respectively  with the same remainder $\nu$, then $aQ_1$ has $j-$capacity $a \xi_1$ and the sum $Q_1 \oplus Q_2$ on $\mathcal{H}_1\oplus \mathcal{H}_2$ has $j-$capacity $(\sqrt[j]{\xi_1}+\sqrt[j]{\xi_2})^j$ both with remainder $\nu$.
	\item \textit{(Independence of restriction)} If $\cV \subseteq \mathcal{H}$ is a subspace of finite codimension then $Q$ has $j-$capacity $\xi$ with remainder $\nu$ if and only if its restriction to $\cV$ has $j-$capacity $\xi$ with remainder $\nu$.
	\item  \textit{(Additivity)} if $Q_1$ has $j-$capacity $\xi$ with remainder $\nu$ and $Q_2$ has $0$ $j-$capacity with remainder of the same order $\nu$, then their sum $Q_1+Q_2$ has the same capacity with remainder $\nu' = \frac{(j+\nu )(j+1)}{j+\nu +1} $ 
\end{itemize}

In the remaining part of this section will be dealing with quadratic forms $Q$ coming from operators of the form  given in \cref{eq: compact part second variation}. 
Suppose that $Z_t$ is a $2n \times k$ matrix which depends piecewise analytically on the parameter $t \in [0,1]$ and  define the following $2n \times 2n $ skew-symmetric matrix:
\begin{equation}
	\label{eq: def J}
	 J = \begin{pmatrix}
		0 &-Id_n \\ Id_n &0
	\end{pmatrix}.
\end{equation}
As $Q$ consider the following quadratic form on $L^2([0,1],\mathbb{R}^k)$:
\begin{equation}
	\label{eq:  def quadratic form}
	Q(v) = \langle v, K(v) \rangle =  \int_0^1 \int_0^t \langle Z_t v(t), J Z_\tau v(\tau)\rangle  d\tau dt.
\end{equation}

\begin{rmrk}
	\label{rmrk: symplectic}
	The operator $K$ and the bilinear form $Q(u,v) = \langle u,K(v) \rangle $ are not symmetric. However the operator:$$K(v) =\int_0^t Z_t^*JZ_\tau v(\tau)d\tau, $$ satisfies \cref{eq: K restricted to V is self adjoint} and becomes symmetric on a finite codimension subspace $\cV$. It is enough to require that the integral $\int_0^1 Z_t v(t) dt$ lies in a Lagrangian subspace of $(\mathbb{R}^{2n},\sigma)$ for any $v \in \cV$. For instance if we consider the fibre (or \emph{vertical} subspace), i.e. the following:
	\begin{equation}
		\label{eq: def fibre}
		\Pi = \{(p,0) : p \in \mathbb{R}^n\} \subset \mathbb{R}^{2n}.
	\end{equation}	
	 Here $\sigma$ denotes the standard symplectic form on $\mathbb R^{2n}$ defined as $\sigma (x,x') = \langle J x,x'\rangle.$ 
\end{rmrk}

Let $f$ be a smooth function on $[0,1]$ and let $k \in \mathbb{N}$, denote by $f^{(k)} = \frac{d^k f}{d t^k}$ the $k-$th derivative with respect to $t$.
For $j\ge 1$ define the following matrix valued functions:
\begin{equation}
	\label{eq: definition Aj}
	A_j(t)= \begin{cases}
		\big(Z_t^{(k)}\big)^*J Z^{(k)}_t \quad &\text{if }j = 2k-1\\
		\big(Z_t^{(k-1)}\big)^*J Z^{(k)}_t \quad &\text{if }j = 2k\\
	\end{cases}
\end{equation}

We use $\rho_t$ to denote any eigenvalue of the matrix $A_{j}(t)$. If $j=2k$, define:
	\begin{equation*}			
	\mu_{t,2k}^+ :=\sum_{\rho_t : \rho_t>0} \sqrt[2k]{ \rho_t } \qquad
	\mu_{t,2k}^- :=\sum_{\rho_t:\rho_t<0} \sqrt[2k]{\vert \rho_t\vert }.
\end{equation*}
For odd indices, $A_{2k-1}$ is skew-symmetric and thus the spectrum is purely imaginary. So we define the function:
\begin{equation*}
	\mu_{t,2k-1} = \sum_{\rho_t : -i\rho_t>0} \sqrt[2k-1]{-i\rho_t}.
\end{equation*}

We are now ready to state the first main result of the section.

\begin{thm}
	\label{thm: eigenvalues of K}
	Let $Q$ be the quadratic form in \cref{eq:  def quadratic form}. $Q$ has either $\infty-$capacity or $j-$capacity with remainder of order $\nu = 1/2$. 
	More precisely, let $j\ge1$ be the lowest integer such that $A_j(t)$ is not identically zero, then
	\begin{itemize}
		\item if $j = 2k-1$, the $(2k-1)-$capacity $\xi$ is given by:
		\begin{equation*}
			\xi = \Bigg(\int_0^1\mu_{t,2k-1} dt \Bigg)^{2k-1},
		\end{equation*}
	and thus for $n\in \mathbb{Z}$ sufficiently large:
	\begin{equation*}
		\lambda_n = \frac{\Big(\int_0^1\mu_{t,2k-1} dt \Big)^{2k-1}}{(\pi n)^{2k-1}} + O(n^{-2k+1/2}).
	\end{equation*}
		\item if $j = 2k$, the $2k-$capacity $(\xi_+,\xi_-)$ is given by:
	\begin{equation*}
		\xi_\pm = \Bigg(\int_0^1\mu^\pm_{t,2k} dt \Bigg)^{2k},
	\end{equation*}
	and thus for $n\in \mathbb{Z}$ sufficiently large:
	\begin{equation*}
		\lambda_n = \frac{\Big(\int_0^1\mu^\pm_{t,2k} dt \Big)^{2k}}{(\pi n)^{2k}} + O(n^{-2k-1/2}). 
	\end{equation*}
	\item if $A_j(t) \equiv 0$ for any $j$ then $Q$ has  $\infty-$capacity. 
	\end{itemize}
	\end{thm}

\begin{rmrk}
	It is worth remarking that in Theorem 1 of \cite{determinant} the order of the remainder for the $1-$capacity was  a little better, $2/3$ and not $1/2$. 
\end{rmrk}

The proof of this result is given in \Cref{section: proof thm 1}. The next theorem gives a characterization of the operators satisfying \cref{eq: K restricted to V is self adjoint,eq: being volterra type} and a geometric interpretation of the $1-$capacity. Before going to the statement let us introduce the following notation. Let $\cA$ denote the skew-symmetric part of $K$:
\begin{equation*}
	\cA = \frac{1}{2}\Big(K-K^*\Big).
\end{equation*}
Let $\Sigma$ be the spectrum of $\cA$ and  $\imm(\cA)$, the image of $\cA$.

\begin{thm}
	\label{thm: characterization of K}
	Let be $K$ an operator satisfying \cref{eq: K restricted to V is self adjoint} and \cref{eq: being volterra type}. Then $\cA$
	has finite rank and completely determines $K$. More precisely, if $\cA$ has rank $2m$ and is represented as:
	\begin{equation*}
		\cA(v)(t):= \frac{1}{2} Z_t^*\cA_0\int_0^1 Z_\tau v(\tau) dt,
	\end{equation*} 
	for a skew-symmetric $2m \times 2m$ matrix  $\cA_0$ and a $2m \times k$ matrix $Z_t$ then:
	\begin{equation}
		\label{eq: factorization of the kernel}
		K(v)(t) = \int_0^tZ_t^*\cA_0Z_\tau v(\tau) d\tau .
	\end{equation}
	Let $\Sigma$ be the spectrum of $\cA$, if the matrix $Z_t$ can be chosen to be piecewise analytic the $1-$capacity of $K$ can be bound by 
	\begin{equation*}
	\xi	 \le 2\sqrt{m}\sqrt{\sum_{\rho \in \Sigma :  -i \rho >0  }-\rho^2} \le2 \sqrt{m}\sum_{\rho \in \Sigma :  -i \rho >0 } \vert\rho \vert.
	\end{equation*}
\end{thm}

		 \section{Proof of \Cref{thm: eigenvalues of K}}
		 \label{section: proof thm 1}
		 
Before going to the proof of \Cref{thm: eigenvalues of K} we still need some auxiliary results. We start with \Cref{lemma taylor Q} to single out the main contributions to the asymptotic of the eigenvalues of $Q$ (the quadratic form defined in \cref{eq:  def quadratic form}).
The first non zero term of the decomposition we give will determine the rate of decaying of the eigenvalues (see \Cref{prop magnitude eigenvalues}).

Before showing this and prove the precise estimates we need to carry out the explicit computation of the asymptotic in some model cases, namely when the matrices $A_j$ are constant. Then we have to show how the $j-$capacity behaves with respect to natural operations such as direct sum of quadratic form or restriction to finite codimension subspaces (\Cref{prop capacity}). 

Let us start with some notation:
\begin{equation*}
	v_k(t) = \int_0^t v_{k-1}(\tau) d\tau, \quad v_0(t) = v(t) \in  L^2([0,1],\mathbb{R}^m)
\end{equation*}
Suppose that the map $t \mapsto Z_t$ is real analytic (or at least regular enough to perform the necessary derivatives) and integrate by parts twice:	
\begin{equation*}
	\begin{split}
		Q (v)&= \int_0^1 \langle Z_t v(t),\int_0^t JZ_\tau v(\tau) d \tau \rangle dt \\ &= \int_0^1 \langle Z_t v(t),JZ_t v_1(t)\rangle - \langle Z_t v(t),\int_0^tJ\dot{Z}_\tau v_1(\tau) d \tau \rangle dt \\
		&= \int_0^1 \langle Z_t v(t),JZ_t v_1(t)\rangle +\langle Z_t v_1(t),J\dot{Z}_t v_1(t)\rangle dt +\\ &\quad \quad +\int_0^1 \langle \dot{Z}_t v_1(t),J\int_0^t\dot{Z}_\tau v_1(\tau) d\tau\rangle dt  - \Big [ \langle \int_0^1 Z_t v(t)dt,J\int_0^1 \dot{Z}_tv_1(t)dt \rangle\Big]
	\end{split}
\end{equation*}

If we impose the condition $\int_0^1 v_t dt =0 \, (\iff v_1(1)=0)$ the term in brackets vanishes:
\begin{equation*}
	 \langle \int_0^1 Z_t v(t)dt,J\int_0^1 \dot{Z}_tv_1(t)dt \rangle =  \langle \int_0^1 Z_t v(t)dt,J Z_1v_1(1) \rangle - \langle \int_0^1 Z_t v(t)dt,J\int_0^1 Z_tv(t)dt \rangle
\end{equation*} and we can write $Q$ as a sum of three terms
\begin{equation*}
	Q(v) = Q_1(v) + Q_2(v)+ R_1(v)
\end{equation*}

In analogy we can make the following definitions: \begin{equation*}
	\begin{split}
		Q_{2k-1}(v) &=  \int_0^1 \langle Z^{(k-1)}_t v_{k-1}(t),JZ^{(k-1)}_t v_{k}(t)\rangle  = \int_0^1 \langle  v_{k-1}(t),A_{2k-1}(t) v_{k}(t)\rangle\\
		Q_{2k}(v) &= \int_0^1\langle Z^{(k-1)}_t v_{k}(t),JZ^{(k)}_t v_{k}(t)\rangle dt = \int_0^1\langle v_{k}(t),A_{2k}(t) v_{k}(t)\rangle dt   \\
		R_k &= \int_0^1 \langle Z^{(k)}_t v_k(t),J\int_0^t{Z}^{(k)}_\tau v_k(\tau) d\tau\rangle dt \\
		V_k &= \{v \in L^2([0,1],\mathbb{R}^m) :  v_l(1)=0, \,\forall \, 0<l \le k\}
	\end{split}
\end{equation*}
Here the matrices $A_j(t)$ are exactly those defined in \cref{eq: definition Aj}.
\begin{lemma}\label{lemma taylor Q}
	For every $j \in \mathbb{N}$, on the subspace $V_j$, the form $Q$ can be represented as\begin{equation}
		\label{eq: taylor Q}
		Q(v) = \sum_{k=1}^{2j} Q_k(v) + R_j(v)
	\end{equation}
	The matrices $A_{2k}(t)$ are symmetric provided  that $\frac{d}{dt} A_{2k-1}(t)\equiv 0$. On the other hand $A_{2k-1}$ is always skew symmetric.
	
	\begin{proof}
		It is sufficient to notice that $R_1(v)$ has the same form as $Q(v)$ but with $v_1$ instead of $v$ and $\dot{Z}_t$ instead of $Z_t$. Thus the same scheme of integration by parts gives the decomposition. 
		
		Notice that $A_{2k}(t) = A^*_{2k}(t)+\frac{d}{dt} A_{2k-1}(t)$ thus the skew-symmetric part of $A_{2k}(t)$ is zero if $A_{2k-1}$ is zero or constant. $A_{2k-1}(t)$ is always skew-symmetric by definition.
	\end{proof}
\end{lemma}

Now we would like to compute explicitly the spectrum of the $Q_j$ when the matrices $A_j$ are constant. Unfortunately describing the spectrum with boundary conditions given by the $V_j$ is quite hard. Already for $Q_4$ the equation determining it cannot be solved explicitly. 

 We will derive the Euler-Lagrange equation for $Q_j$ and turn instead to periodic boundary conditions for which everything becomes very explicit and show how to relate the solution for the two boundary value problems we are considering. Let us write down the Euler-Lagrange equations for the forms $Q_j$. If $j = 2k$ integration by parts yields:
 \begin{equation*}
 	\begin{split}
 		Q_{2k}(v)- \lambda \vert\vert v\vert\vert^2=& \int_0^1 \langle v_k(t), A_{2k} v_k(t)\rangle- \lambda \langle v_0(t),v_0(t) \rangle dt \\
 		=& \int_0^1 \langle v_0(t),(-1)^kA_{2k}v_{2k}(t)-\lambda v_0(t)\rangle dt+ \\& \quad +\sum_{r=0}^{k-1}(-1)^{r} \Big[\langle v_{k-r}(t),A_{2k}v_{k+r+1}(t)\rangle\Big]_0^1
 	\end{split}
 \end{equation*} 
 
Notice that the boundary terms vanish identically if we impose the vanishing of $v_j$ for $1 \le j\le k$ at boundary points.
 
 We change notation and define  $w(t) = v_{2k}(t)$ and $w^{(j)} (t) = \frac{d^j}{dt^j}(w(t))$. The new equations are:
 \begin{equation*}
 	w^{(2k)}(t) = \frac{(-1)^{k}}{\lambda} A_{2k} w(t)
 \end{equation*}
 
 We can perform a linear change of coordinates that diagonalizes $A_{2k}$ to reduce to $m$ $1-$dimensional systems. Imposing periodic boundary conditions, we are thus left with the following boundary value problem:
 \begin{equation}
 	\label{eq: boundary value problem even}
 	w^{(2k)}(t) = \frac{(-1)^{k}\mu}{\lambda} w(t) \quad w^{(j)}(0)=w^{(j)}(1) \text{ for } 0\le j\le 2k-1
  \end{equation}
 
 The case of odd $j$ is very similar, in fact $Q_{2k-1}(v)$ can be rewritten as:
 \begin{equation*}
 	\begin{split}
 		Q_{2k-1}(v)-\lambda \vert\vert v\vert\vert^2 = & \int_0^1 \langle v_{k-1}(t), A_{2k-1} v_k(t)\rangle- \lambda \langle v_0(t),v_0(t) \rangle dt \\
 		=& \int_0^1 \langle v_0(t), (-1)^{k-1}A_{2k-1}v_{2k-1}(t)-\lambda v_0\rangle dt + \textit{ b.t.}
 	\end{split}
 \end{equation*}
 
 Here by $b.t.$ we mean boundary terms as the one appearing in the previous equation. They again disappear if we assume that $v_{j}\in V_j$.
 Thus we end up with a boundary value problem similar to the one we had before with the difference that now the matrix $A_{2k-1}$ is skew-symmetric.
 \begin{equation*}
 	w^{(2k-1)}(t) = \frac{(-1)^{k-1}}{\lambda} A_{2k-1} w(t)
 \end{equation*}
 
 If we split the space into the kernel and invariant subspaces on which $A_{2k-1}$ is non degenerate we can decompose $Q_{2k-1}$ as a direct sum of two-dimensional forms. Imposing periodic boundary conditions, we end up with the following boundary value problems:
 \begin{equation}
 	\label{eq: boundary value odd}
 	\begin{cases}
 		w_1^{(2k-1)}(t) &=- \frac{(-1)^{(k-1)}\mu}{\lambda}w_2\\
 		w_2^{(2k-1)}(t) &= \frac{(-1)^{(k-1)}\mu}{\lambda} w_1
 	\end{cases} \quad \begin{cases}
 w_1^{(j)}(0) = w_1^{(j)}(1), \\  w_2^{(j)}(0) = w_2^{(j)}(1)
\end{cases} \text{ for } 0\le j\le 2k-2.
 \end{equation}

\begin{lemma}
	\label{lemma: boundary value problem}
	The boundary value problem in \cref{eq: boundary value problem even} has a solution if and only if $$\lambda \in \Big\{\frac{\mu}{(2\pi r)^{2k}}: r \in \mathbb{N}\Big\}.$$
	Moreover any such $\lambda$ has multiplicity $2$. In particular, the decreasing sequence of $\lambda$ for which \cref{eq: boundary value problem even} has  solutions satisfies:
	\begin{equation*}
		\lambda_r = \frac{\mu}{(2\pi \lceil r/2 \rceil)^{2k}} = \frac{\mu}{(\pi r)^{2k}} + O(r^{-(2k+1)}), \quad r \in \mathbb{N} 
	\end{equation*}
	
	Similarly the boundary value problem in \eqref{eq: boundary value odd} has a solution if and only if:
	\begin{equation*}
		\lambda \in \Big\{\frac{\vert\mu \vert}{(2\pi r)^{2k-1}}: r \in \mathbb{Z}\Big\}
	\end{equation*}
	and any such $\lambda$ has again  multiplicity $2$. The  monotone rearrangement of $\lambda$ for which there exists  a solution to the boundary value problem is:
	\begin{equation*}
		\lambda_r = \frac{\vert\mu\vert }{(2\pi \lceil r/2 \rceil)^{2k-1}} = \frac{\vert \mu \vert}{(\pi r)^{2k-1}} + O(r^{-(2k)}), \quad r \in \mathbb{Z} 
	\end{equation*}
	
	\begin{proof}
	Any solution of the equation $w^{(2k)}(t) = \frac{(-1)^{k}\mu}{\lambda} w(t)$ can be expressed as a combination of trigonometric and hyperbolic functions with the appropriate frequencies.
	
	Without loss of generality we can assume $\mu>0$, we have to consider two separate cases: 
	\caso{$k$ even and $\lambda>0$ or $k $ odd and $\lambda<0$ }
	In this case the quantity $(-1)^k\mu \lambda^{-1}>0$. If we define $a^{2k} =(-1)^k\mu \lambda^{-1}>0$ for $a>0$, we have to solve:
	\begin{equation}
		\label{eq: linear boundary value problem eigenvalues}
		w^{(2k)} (t)= a^{2k} w(t), \qquad w^{(j)}(0)=w^{(j)}(1),\,\, 0\le j<2k.
	\end{equation}

	A base for the space of solutions to the $ODE$ is then $\{e^{\omega^j a t}:  \omega = e^{i\pi/k}\}$. For us it will be more convenient to switch to a real representation of the space of solutions. Notice the following symmetry of the even roots of $1$, if $\eta$ is a root of $1$ different form $\pm 1, \pm i$ then $\{\eta, \bar \eta,-\eta, -\bar\eta\}$ are still distinct roots of $1$ (this is also a Hamiltonian feature of the problem).
		
	If we write $\eta = \eta_1+i \eta_2$, this symmetry implies that the space generated by $\{e^{\eta t }, e^{\bar \eta t },e^{-\eta t },e^{-\bar \eta t}\}$ is the same as the space generated by $$\{\sin(\eta_2 t)\sinh(\eta_1 t),\sin(\eta_2 t)\cosh(\eta_1 t),\cos(\eta_2 t)\sinh(\eta_1 t),\cos(\eta_2 t)\cosh(\eta_1 t)\}.$$
	
	Let us rescale these functions by $a$ (so that they solve \cref{eq: linear boundary value problem eigenvalues}) and call their linear span $U_\eta$, we then define $U_1$ to be the span of $\{\sinh(t),\cosh(t)\}$ and $U_i = \{\sin(t),\cos(t)\}$. Note that $U_i$ appears if and only if $k$ is even. 
		
		Thus the solution space for our problem is the space $\bigoplus_\eta U_\eta$ where $\eta$ ranges over the set $E = \{\eta : \Re(\eta)\ge0, \Im(\eta)\ge0,\eta^{2k}=1 \}$.  
		
		Now we have to impose the boundary conditions. Notice that, if $k$ is even then $U_i$ is made of periodic functions, so they are always solutions. We can look for more on the complement $\bigoplus_{\eta \ne i} U_\eta$. Suppose by contradiction that $w$ is one of such solutions. Write $w = \sum_\eta w_\eta$ with $w_\eta \in U_\eta$ and let $b$ be the $\sup \{\Re(\eta) : \eta \in E, w_\eta \ne 0\}$. It follows that either $\sinh(b\, a t)$ or $\cosh(b \, a t)$ is present in the decomposition of $w$. It follows that:
		$$w(t) = \sinh(b \, a t) \frac{w(t)}{\sinh( b \, a t )} = \sinh(b \, a t) g(t), \quad  0 \not \equiv \vert g(t) \vert < C \text{ for $t$ large enough}$$
		and so $\vert w \vert$ is unbounded as $t \to + \infty$ (or $-\infty$) and thus $w$ is not periodic. It follows that there are periodic solutions only if  $k$ is even (and thus $\lambda>0$) and $a = 2\pi r = \sqrt[2k]{\frac{\mu}{\lambda}}$. Notice that we have two independent solutions, so if we order the solution decreasingly we have:
		\begin{equation*}
			\lambda_r = \frac{\mu}{(2\pi \lceil r/2 \rceil)^{2k}}, \quad r \in \mathbb{N}
		\end{equation*}
		\caso{$k$ odd and $\lambda>0$ or $k$ even and $\lambda<0$}
		In this case we have to look at the roots of $-1$ but the argument is very similar. If $k$ is even there are no solutions, since you lack purely imaginary frequencies.  If $k$ is odd, set $\vert\mu\lambda^{-1}\vert=a^{2k}$, then the boundary value problem is:
		\begin{equation*}		
			w^{(2k)} (t)=- a^{2k} w(t) \qquad w^{(j)}(0)=w^{(j)}(1),\, 0\le j<2k.
		\end{equation*}
		
		The roots of $-1$ are just the roots of $1$ rotated by $i$. Now the space of solutions is $\bigoplus_{\eta\ne 1} U_\eta $. We find again two independent solutions, if we order them we get:
		\begin{equation*}
			\lambda_{r} = \frac{\mu}{(2\pi \lceil r/2 \rceil)^{2k}}, \quad r \in \mathbb{N}
		\end{equation*}

		Notice that positive $\mu$ gives rise to positive solutions. Thus if we consider $\mu <0$, we get the same result but with switched signs.

		We can reduce the odd case (\cref{eq: boundary value odd}) to the even one. Consider the $1-$dimensional equation of twice the order, i.e.:
		\begin{equation*}
			w_1^{2(2k-1)}(t) =- \frac{\mu^2}{\lambda^2} w_1
		\end{equation*}
		Now, the discussion above tells us that there are exactly two independent solutions with periodic boundary conditions whenever $\lambda$ satisfies $\sqrt[2k-1]{\frac{\mu}{\vert\lambda\vert}} = 2 r \pi$. It follows that again there are two independent solutions, this times for both signs of $\lambda$. If we order them we get:
		\begin{equation*}
			\lambda_r = \frac{\mu}{(2 \pi \lceil r/2 \rceil)^{2k-1}}, \quad \lambda_{-r} = \frac{\mu}{(2 \pi \lfloor r/2 \rfloor)^{2k-1}}, \quad r \in \mathbb{N}
		\end{equation*}
	\end{proof}
\end{lemma}

\begin{prop}
	\label{prop: counting functions maslov index}
	Let $\mu>0$ and $s \in (0,+\infty)$, denote by $\eta_s$ the number of solutions of \cref{eq: boundary value problem even} with $\lambda$ greater than $s$ and similarly denote by $\omega_s$ be the number of solutions with $\lambda$ bigger than $s$ of:
	\begin{equation}
		\label{eq: boundary value problem right bc}
		w^{(2k)}(t) = \frac{(-1)^{k}\mu}{\lambda} w(t), \quad w^{(j)}(0) = w^{(j)}(1)=0,\quad k\le j\le 2k-1
	\end{equation}
	Then $\vert \omega_s-\eta_s \vert \le 2k$. The same conclusion holds for \cref{eq: boundary value odd}.
	\begin{proof}
		The result follows from standard results about Maslov index of a path in the Lagrange Grassmannian. References on the topic can be found in \cite{beschastnyi_morse, beschastnyi_1d, agrachev_quadratic_paper}. Let us illustrate briefly the construction. Let $(\Sigma, \sigma)$ be a symplectic space, the Lagrange Grassmannian is the collection of Lagrangian subspaces of $\Sigma$ and it has a structure of smooth manifold. For any Lagrangian subspace $L_0$ we define  the \emph{train} of $L_0$ to be the set: $T_{L_0}=\{L \text{ Lagrangian}: L \cap L_0 \ne (0)\}$. $T_{L_0}$ is a stratified set, the biggest stratum has codimension $1$ and is endowed with a co-orientation. If $\gamma$ is a smooth curve with values in the Lagrangian Grassmannian (i.e. a smooth family of Lagrangian subspaces) which intersects transversally $T_{L_0}$ in its smooth part, one defines an intersection number by counting the intersection points weighted with a plus or minus sign depending on the co-orientation. 
		Tangent vectors at a point $L$ of the Lagrange Grassmannian (which is a subspace of $\Sigma$) are naturally interpreted as quadratic forms on $L$. We say that a curve is \emph{monotone} if at any point its velocity is either a non negative or a non positive quadratic form.  		
		For monotone curves, Maslov index counts the number of intersections with the train up to sign. For generic continuous curves it is defined via a homotopy argument.   
		
		Denote by $\mathrm{Mi}_{L_0}(\gamma)$ the Maslov index of a curve $\gamma$ and $L_1$ be another Lagrangian subspace. In \cite{agrachev_quadratic_paper} the following inequality is proved:
		\begin{equation}
			\label{eq: inequality maslov}
			\vert \mathrm{Mi}_{L_0}(\gamma)-\mathrm{Mi}_{L_1}(\gamma)\vert \le \frac{\dim(\Sigma)}{2}
		\end{equation}
	
		Let us apply this results to our problem. First of all let us produce a curve in the Lagrange Grassmannian whose Maslov index coincides with the counting functions $\omega_s$ and $\eta_s$. The right candidate is the graph of the fundamental solution of $w^{(2k)}(t) = \frac{(-1)^{k}\mu}{\lambda} w(t)$. 
		
		We write down a first order system on $\mathbb{R}^{2k}$ equivalent to our boundary value problem, if we call the coordinates on $\mathbb{R}^{2k}$ $x_j$, set:
		\begin{equation*}
			x_{j+1} (t)= w^{(j)}(t) \Rightarrow \dot{x}_j = x_{j+1} \text{ for }1\le j\le 2k-1, \quad \dot x_{2k} = \frac{(-1)^k \mu}{\lambda}x_1.
		\end{equation*}
		For simplicity call $\frac{(-1)^k \mu}{\lambda} = a$, the matrix we obtain has the following structure:
		\begin{equation*}
			A_\lambda= \begin{pmatrix}
				0 &       &		&	a\\
				1 &0	  &		&	\\
				 &\ddots &\ddots&	\\
				 &		  &1		&0	
			\end{pmatrix}
		\end{equation*}
	This matrix is not Hamiltonian with respect to the standard symplectic form on $\mathbb{R}^{2k}$ but is straightforward to compute a similarity transformation that sends it to an Hamiltonian one (recall that we already used that $A_\lambda$ has the spectrum of an Hamiltonian matrix). Moreover the change of coordinates can be chosen to be block diagonal and thus preserves the subspace $B = \{x_j = 0,  k\le j\}$, which remains Lagrangian too. Since later on we will have to show that the curve we consider is monotone we will give this change of coordinates explicitly. Define the matrix $S$ setting $S_{i,k-i+1}= (-1)^{i-1}$ and zero otherwise. It is a matrix that has alternating $\pm 1$ on the anti-diagonal. Define the following $2k\times 2k$ matrices:
	\begin{equation*}
		G = \begin{pmatrix}
			1 &0\\0 &S
		\end{pmatrix} \quad G^{-1} = \begin{pmatrix}
		1 &0 \\ 0 & (-1)^kS
	\end{pmatrix} \quad \hat{A}_\lambda = G A_\lambda G^{-1}
	\end{equation*} 
	Set $N$ to be the lower triangular $k\times k $ shift matrix (i.e. the left upper block  of $A_\lambda$ above) and $E$ the matrix with just a $1$ in position $(1,k)$ (i.e. the left lower block of $A_\lambda$). The new matrix of coefficients is:
	\begin{equation*}
		\hat A_\lambda =\begin{pmatrix}
			N &a(-1)^k ES \\ SE &-N^*
		\end{pmatrix} \quad ES = \mathrm{diag}(0,\dots,0,1), \quad  SE = \mathrm{diag}(1,0,\dots,0).
	\end{equation*}
	Now we are ready to define our curve. First of all the symplectic space we are going to use is $(\mathbb{R}^{4k}, \sigma\oplus(- \sigma))$ where $\sigma$ is the standard symplectic form, in this way graphs of symplectic transformation are Lagrangian subspaces. Sometimes we will denote the direct sum of the two symplectic forms with opposite signs with $\sigma \ominus \sigma$ too. Let $\Phi_\lambda$ be the fundamental solution of $\dot \Phi_\lambda^t = \hat A_\lambda \Phi_\lambda^t$ at time $t =1$. Consider its graph:
	\begin{equation*}
		\gamma: \lambda \mapsto \Gamma(\Phi^1_\lambda)= \Gamma(\Phi_\lambda), \quad  \lambda \in (0,+\infty)
	\end{equation*} 
	
	Once we prove that $\gamma$ is monotone, is straightforward to check that $\mathrm{Mi}_{B \times B} (\gamma\vert_{[s,+\infty)})$ counts the number of solutions to boundary value problem given in \cref{eq: boundary value problem right bc} for $\lambda \ge s$ and similarly $\mathrm{Mi}_{\Gamma(I)} (\gamma\vert_{[s,+\infty)})$ counts the solutions of \cref{eq: boundary value problem even} for $\lambda\ge s$. Here $\Gamma(I)$ stands for the graph of the identity map (i.e. the diagonal subspace).
	
	Let us check that the curve is monotone. As already mentioned, tangent vectors in the Lagrange Grassmannian can be interpreted as quadratic forms. Being monotone means that the following quadratic form is either non negative or non positive:
	\begin{equation*}
		\big(\partial_\lambda\gamma\big)(\xi)= \sigma (\Phi_\lambda \xi,\partial_\lambda \Phi_\lambda \xi), \quad \xi \in \mathbb{R}^{2k}
	\end{equation*} 
	We use the ODE for $\Phi_\lambda(t)$ to prove monotonicity:
	\begin{equation*}
		\begin{split}
		 \sigma (\Phi_\lambda \xi,\partial_\lambda \Phi_\lambda \xi) &=\int_0^1 \frac{d}{dt}\big(\sigma (\Phi^t_\lambda \xi,\partial_\lambda \Phi^t_\lambda \xi)\big)dt + \sigma (\Phi^0_\lambda \xi,\partial_\lambda \Phi^0_\lambda \xi)\\ &=\int_0^1\sigma(\hat A_\lambda \Phi^t_\lambda\xi,\partial_\lambda \Phi^t_\lambda \xi) + \sigma( \Phi^t_\lambda\xi,\big(\partial_\lambda\hat A_\lambda \, \Phi^t_\lambda +\hat{A}_\lambda \partial_\lambda \Phi^t_\lambda \big)\xi)dt \\
		 &= \int_0^1\sigma( \Phi^t_\lambda\xi,\partial_\lambda\hat A_\lambda \, \Phi^t_\lambda \xi)dt
	 \end{split}
	\end{equation*}
	
	Where we used the facts that $\partial_\lambda \Phi^0_\lambda = \partial_\lambda Id =0$  and that $\hat A_\lambda$ is Hamiltonian and thus $J\hat A_\lambda = -\hat A_\lambda^* J $ to cancel the first and third term. It remains to check $J\partial_\lambda \hat{A}_\lambda$. It is straightforward to see that it is a diagonal matrix with just a non zero entry, thus is either non negative or non positive. So $\partial_\lambda \gamma$ is either non positive or non negative being the integral of a non positive or non negative quantity (the sign is independent of $\xi$).
	
	Now the statement follows from inequality \eqref{eq: inequality maslov}.
	\end{proof}
\end{prop}

We are finally ready to compute the asymptotic for $Q_j$ when the matrix $A_j$ is constant. The next Proposition translate the estimate on the counting functions $\eta_s$ and $\omega_s$ defined in \Cref{prop: counting functions maslov index} to an estimate for the eigenvalues.

\begin{prop}
	Let $Q_j$ be any of the forms appearing in \cref{eq: taylor Q}.
	\begin{itemize}
		\item Suppose $j =2k$ and $Q_{2k}(v) = \int_0^1 \langle A_{2k}v_k, v_k \rangle dt$ with $A_{2k}$ symmetric and constant and let $\Sigma_{2k}$ be its spectrum. Define $$\xi_+ = \left(\sum_{\mu \in \Sigma_{2k},\mu>0} \sqrt[j]{\mu} \right)^j \text{ and } \xi_- =\left(\sum_{\mu \in \Sigma_{2k},\mu<0} \sqrt[j]{\vert\mu\vert} \right)^j. $$ Then $Q_{2k}$ has capacity $(\xi_+,\xi_-)$ with remainder of order one. 
		Moreover, if $A_{2k}$ is $m\times m$ and $r \in \mathbb{N}$, for $r \ge m k$
		\begin{equation}
			\label{eq: bound eigenvalues linear even}						\frac{\xi_+}{\pi^j (r-2mk-p(r))^j} \ge \lambda_{r} \ge\frac{\xi_+}{\pi^j(r+2mk+p(r))^j}  
		\end{equation}
		where $p(r) = 0$ if $r$ is even or $p(r)=1$ if $r$ is odd. Similarly for negative $r$ with $\xi_-$.
			\item Suppose $j =2k+1$ and $Q_{2k+1}(v) = \int_0^1 \langle A_{2k+1}v_{k-1}, v_{k} \rangle dt$ with $A_{2k+1}$ skew-symmetric and constant and let $\Sigma_{2k+1}$ be its spectrum. Define $$\xi = \left(\sum_{\mu \in \Sigma_{2k+1},-i\mu>0} \sqrt[j]{-i\mu} \right)^j.$$ Then $Q_{2k+1}$ has capacity $\xi$ with remainder of order one. 
			Moreover , if $A_{2k}$ is $m\times m$ and $r \in \mathbb{Z}$, for $\vert r\vert \ge m k$
			\begin{equation}
				\label{eq: bound eigenvalues linear odd}						\frac{\xi}{\pi^j (r-2mk-p(r))^j} \ge \lambda_{r} \ge\frac{\xi}{\pi^j(r+2mk+p(r))^j} . 
			\end{equation}
	\end{itemize}
	\begin{proof}
		First of all we consider $1-$dimensional system and we write the inequality $\vert\eta_s-\omega_s\vert$ as an inequality for the eigenvalues. Notice that if we have two integer valued function $f,g : \mathbb{R}\to \mathbb{N}$ and an inequality of the form:
		\begin{equation*}
			g(s)\ge \#\{\lambda \text{ solutions of \cref{eq: boundary value problem right bc} }: \lambda \ge s \} \ge f(s),
		\end{equation*}
		it means that we have at least $f(s)$ solutions bigger than $s$ and at most $g(s)$. This implies that the sequence of ordered eigenvalues satisfies:\begin{equation*}
			\lambda_{f(s)} \ge s, \quad \lambda_{g(s)}\le s.		
			\end{equation*}
		Now we compute this quantities explicitly. In virtue of \Cref{prop: counting functions maslov index} we can take as upper/lower bounds for the counting function $g(s) = \eta_s + 2k$ and $f(s) = \eta_s-2k$. We choose the point $s = \frac{\mu}{(2\pi r)^j}$. It is straightforward to see that:
		\begin{equation*}
			\eta_s\Big\vert_{s =\frac{\mu}{(2\pi r)^j}} = 2 \#\{l \in \mathbb{N} : \frac{\mu}{(2\pi l)^j}\ge\frac{\mu}{(2\pi r)^j}\} = 2 r.
		\end{equation*}
		And thus we obtain:
		\begin{equation*}
			\lambda_{2(r-k)} \ge \frac{\mu}{(2\pi r)^j}, \quad \lambda_{2(r+k)}\le \frac{\mu}{(2\pi r)^j}.		
		\end{equation*}
	Now if we change the labelling we find that , for $l\ge k$:
	\begin{equation*}
		\frac{\mu}{(2\pi (l-k))^j}\ge \lambda_{2l} \ge \frac{\mu}{(2\pi (l+k))^j}.		
	\end{equation*} By definition $\lambda_{2l}\ge \lambda_{2l+1}\ge \lambda_{2l+2}$ and thus we have a bound for any index $r \in \mathbb{N}$.

	Now we consider $m-$dimensional system, notice that we reduced the problem, via diagonalization, to the sum of $m$ $1-$dimensional systems. Thus our form $Q_j$ is always a direct sum of $1-$ dimensional objects. We show now how to recover the desired estimate for the sum of quadratic forms.
	
	First of all observe that counting functions are additive with respect to direct sum. In fact, if $Q =  \oplus_{i=1}^m Q_i$, $\lambda$ is an eigenvalue of $Q$ if and only if it is an eigenvalue of $Q_i$ for some $i$. We proceed as we did before. Suppose that $Q_a$ is $1-$dimensional and $Q_a (v) = \int_0^1 \mu_a \vert v_k(t)\vert ^2 dt $. Let us compute $\eta_s$ in the point $s_0 = (\sum_{i=1}^{m} \sqrt[j]{\mu_i})^j/(2\pi l)^j$:
	\begin{equation*}
		2\#\left\{r\in\mathbb{N} : \frac{\mu_a}{(2\pi r)^j} \ge\frac {(\sum_{i=1}^{m} \sqrt[j]{\mu_i})^j}{(2\pi l)^j} \right\} = 2\#\left\{r\in\mathbb{N} : \frac{\sqrt[j]{\mu_a}}{(\sum_{i=1}^{m} \sqrt[j]{\mu_i})r } \ge\frac {1}{l} \right\}
	\end{equation*}
	Set for simplicity $c_a = \frac{\sqrt[j]{\mu_a}}{(\sum_{i=1}^{m} \sqrt[j]{\mu_i})}$, it is straightforward to see that the cardinality of the above set is $\# \{r \in \mathbb{N}:r\le c_a l\} = \lfloor c_a l \rfloor $. Now we are ready to prove the estimates for the direct sum of forms. Adding everything we have:
	\begin{equation*}
		2\sum_{a=1}^m (\lfloor c_a l \rfloor +k) \ge \#\Big\{ \text{eigenvalues of } Q \ge \frac {(\sum_{i=1}^{m} \sqrt[j]{\mu_i})^j}{(2\pi l)^j} \Big\} = 2\sum_{a=1}^{m}(\lfloor c_a l \rfloor -k)
	\end{equation*}
	It is clear that $\sum_{a=1}^m c_a =1$ and that $l+mk \ge \sum_{a=1}^m (\lfloor c_a l \rfloor +k)$, similarly $\sum_{a=1}^m (\lfloor c_a l \rfloor +k) \ge l-m(k+1)$ since $\lfloor c_a l \rfloor\ge c_al-1$. Rewriting for the eigenvalues with $l\ge mk$ we obtain:
	\begin{equation*}
		\frac{(\sum_{i=1}^{m} \sqrt[j]{\mu_i})^j}{(2\pi(l-mk))^j} \ge \lambda_{2l} \ge\frac{(\sum_{i=1}^{m} \sqrt[j]{\mu_i})^j}{(2\pi(l+mk))^j}.
	\end{equation*}
It is straightforward to compute the bounds in \cref{eq: bound eigenvalues linear even,eq: bound eigenvalues linear odd} observing again $\lambda_{2l}\ge \lambda_{2l+1}\ge \lambda_{2l+2}$.
	\end{proof}
\end{prop}

\begin{rmrk}
	The shift $m$ appearing in \cref{eq: bound eigenvalues linear even,eq: bound eigenvalues linear odd} is due to the fact we are considering the direct sum of $m$ quadratic forms. It is worth noticing that this does not depend on the fact that we are considering a quadratic form on $L^2([0,1],\mathbb
	R^m)$ and the estimates in \cref{eq: bound eigenvalues linear even,eq: bound eigenvalues linear odd} hold whenever we consider the direct sum of $m$ $1-$dimensional forms with constant coefficients. This consideration will be used in the proof of \Cref{thm: eigenvalues of K} below.
\end{rmrk}
Now we prove some properties of the capacities which are closely related to the explicit estimate we have just proved for the linear case. As done so far we state the proposition for ordered positive eigenvalues. An analogous statement is true for the negative ones.
\begin{prop}
	\label{prop capacity}
	Suppose that $Q$ is a quadratic form on an Hilbert space and let $\{\lambda_n\}_{n \in \mathbb{N}}$ be its positive ordered eigenvalues. Suppose that:
	\begin{equation*}
		\lambda_n = \frac{\zeta}{n^j} + O(n^{-j-\nu}) \quad \nu >0, j \in \mathbb{N} \text{ as } n \to +\infty.
	\end{equation*}
	\begin{enumerate}
		\item Then for any such $Q_i$ on a Hilbert space $\mathcal{H}_i$ the direct sum $Q = \oplus_{i=1}^m Q_i$ satisfies:
		\begin{equation*}
			\lambda_n = \Big(\sum_{i=1}^{m}\frac{\sqrt[j]{\zeta_i}}{n}\Big)^j + O(n^{-j-\nu}) \quad \nu >0, j \in \mathbb{N} \text{ as } n \to +\infty.
		\end{equation*} 
		\item  Suppose that $U$ is a subspace of codimension $d<\infty$ then \begin{equation*}
			\lambda_n(Q\vert_U) = \frac{\zeta}{n^j} + O(n^{-j-\nu}) \iff\lambda_n(Q) = \frac{\zeta}{n^j} + O(n^{-j-\nu}),
		\end{equation*}
		as $ n \to +\infty$.
		\item Suppose that $Q$ and $\hat{Q}$ are two quadratic forms. Suppose that $Q$ is as at the beginning of the proposition and $\hat{Q}$ satisfies:
		\begin{equation*}
			\lambda_n(\hat{Q}) = O(n^{j+\mu}) \quad \mu >0,  \text{ as } n \to + \infty.
		\end{equation*}
		Then the sum $Q' = Q+\hat{Q}$ satisfies:
		\begin{equation*}
			\lambda_n(Q') = \frac{\zeta}{n^j} + O(n^{j+\nu'}), \quad \nu' = \min\{\frac{j+\mu}{j+\mu+1}(j+1),j+\nu\}.
		\end{equation*}
	\end{enumerate}
	\begin{proof}
		The asymptotic relation can be written in terms of a counting function. Take the $j-$th root of the eigenvalues of $Q_i$, then it holds that \begin{equation*}
			\# \{n \in \mathbb{N}\, \vert \, 0\le \frac{1}{\sqrt[j]{\lambda_n}}\le k\} = \sqrt[j]{\zeta_i} k+ O(k^{1-\nu})
		\end{equation*}
		So summing up all the contribution we get the estimate in $i)$.
		
		The min-max principle implies that we can control the $n-$th eigenvalue of $Q\vert_U$ with the $n-$th and $(n+d)$-th eigenvalue of $Q$  i.e.:
		\begin{equation*}
			\lambda_{n}(Q\vert_U) \le\lambda_n(Q)\le \lambda_{n-d}(Q\vert_U) \le\lambda_{n-d}(Q) 
		\end{equation*}
		So, if the codimension is fixed, it is equivalent to provide and estimate for the eigenvalues $Q$ or for those of $Q\vert_U$.

		For the last point we use Weyl law. We can estimate the $i+j$-th eigenvalue of a sum  of quadratic forms with the sum of the $i-$th and the $j$-th eigenvalues of the summands. Write, as in \cite{determinant}, $Q'$ as $Q$+$\hat{Q}$ and $Q$ as $Q'$+$(-\hat{Q})$. and choose $i = n-\lfloor n^\delta \rfloor $ and $j = \lfloor n^\delta \rfloor$ in the first case and $i = n$ and $j = \lfloor n^\delta \rfloor$ in the second. This implies:
		\begin{equation*}
			\lambda_{n+\lfloor n^\delta \rfloor}(Q)+\lambda_{\lfloor n^\delta \rfloor}(\hat{Q}) \le \lambda_n(Q')\le \lambda_{n-\lfloor n^\delta \rfloor}(Q)+\lambda_{\lfloor n^\delta \rfloor}(\hat{Q})
		\end{equation*}
		The best remainder is computed as $\nu' = \max_{\delta\in (0,1)}\min\{(j+\mu)\delta,j+1-\delta,j+\nu\}$. 
	\end{proof}
\end{prop}
 
Collecting all the facts above we have the following estimate on the decaying of the eigenvalues of $Q_j$, independently of any analyticity assumption of the kernel.

\begin{prop} 
	\label{prop magnitude eigenvalues}
	Take $Q_j$ as in the decomposition of lemma \eqref{lemma taylor Q}. Then the eigenvalues of $Q_j$ satisfy:
	\begin{equation*}
		\lambda_n(Q_j) = O\Big(\frac{1}{n^j} \Big) \quad \text{ as } n \to \pm \infty 
	\end{equation*}
	Moreover for any $k \in \mathbb{N}$ and for any $0 \le s \le k$ the forms $Q_{2k+1}$ and $Q_{2k}$ have the same first term asymptotic as the forms:
	\begin{equation*}
		\begin{aligned}
			\hat{Q}_{2k+1,s}(v) = (-1)^{s}\int_0^1 \langle A_{2k+1} v_{k+1+s}(t),v_{k-s}(t) \rangle dt \\
			\hat{Q}_{2k,s}(v) = (-1)^{s}\int_0^1 \langle A_{2k} v_{k+s}(t),v_{k-s}(t) \rangle dt
		\end{aligned}
	\end{equation*}
	\begin{proof}
		Let's start with even case, $j=2k$. It holds that:
		\begin{equation*}
			\vert Q_{2k}(v)\vert = \vert\int_0^1\langle A_t v_k(t),v_k(t)dt \vert \le C \int_0^1 \langle v_k(t),v_k(t) \rangle dt
		\end{equation*}
		Where $C = \max_t \vert\vert A_t\vert\vert$. By comparison with the constant coefficient case  we get the bound.
		
		Suppose now that $j = 2k-1$. As before there is a constant $C$ such that
		\begin{equation*}
			\vert Q_{2k}(v)\vert = \vert\int_0^1\langle A_t v_k(t),v_{k+1}(t)dt \vert \le C \Vert v_k \Vert_2 \Vert v_ {k+1}\Vert_2 
		\end{equation*}
		
		Consider now the following quadratic forms on $L^2([0,1],\mathbb{R}^k)$: $$F_k(v) = \int_0^1\vert \vert v_k(t)\vert\vert^2 dt=\Vert v_k \Vert_2  ^2, \quad F_{k+1}(v) = \int_0^1\vert \vert v_{k+1}(t)\vert\vert^2 dt= \Vert v_ {k+1}\Vert_2 ^2$$
		Define $V_n = \{v_1, \dots, v_n\}^{\perp}$ where $v_i$  are linearly independent eigenvectors of $F_k$ associated to the first $n$ eigenvalues $\lambda_1\ge \dots \ge\lambda_n$. Similarly  define $U_n = \{u_1, \dots, u_n\}^{\perp} $ to be  the orthogonal complement to the eigenspace associated to the first $n$ eigenvalues of $F_{k+1}$.
		It follows that:
		\begin{equation*}
			\lambda_{2n}(Q_{2k+1})\le \max_{v \in V_n\cap U_n} C \Vert v_{k}\Vert_2\Vert v_{k+1}\Vert_2 \le  C \max_{v \in V_n}  \Vert v_{k}\Vert_2\max_{v \in U_n}  \Vert v_{k+1}\Vert_2
		\end{equation*}
		
		We already have an estimate for the eigenvalues of $F_k$ and $F_{k+1}$ since we have already dealt with constant coefficients case. In virtue of the choice of the subspace $V_n$ and $U_n$, the maxima in the right hand side are the square roots of the $n-th$ eigenvalues of the respective forms. Thus one gives a contribution of order $n^{-k}$ and the other of order $n^{-k-1}$ and the first part of the proposition is proved.
		
		For the second part, without loss of generality suppose that $j=2k$. The other case is completely analogous. 
		\begin{equation*}
			\begin{split}
				Q_{2k}(v) &= \int_0^1\langle  v_{k},A_t v_{k} \rangle dt = \int_0^1 \langle v_k, \int_0^t A_\tau v_{k-1}(\tau)+\dot{A}_\tau v_{k}(\tau) d\tau \rangle dt\\
				&= - \int_0^1 \langle v_{k+1}(t),A_tv_{k-1}(t)+\int_0^1\langle v_{k+1}(t),\dot{A}_t v_{k}(t)\rangle dt \\
			\end{split}
		\end{equation*}
		The second term above is of higher order by the first part of the lemma and so iterating the integration by parts on the first term at step $s$ we get that:
		\begin{equation*}
			\begin{aligned}
				\int_0^1 \langle v_{k+s}(t),A_t v_{k-s}(t) \rangle dt = - \int_0^1 \langle v_{k+s+1}(t),A_t v_{k-s-1}(t) \rangle dt + \\+ \int_0^1\langle v_{k+s+1}(t), \dot{A}_\tau v_{k-s}(t)  \rangle dt
			\end{aligned}
		\end{equation*}
		
		The second term of the right hand side is again of order $n^{2k+1}$, this can be checked in the same way as in the first part of the proposition. This finishes the proof.
	\end{proof}
\end{prop}

Now we prove the main result of this section:

\begin{proof}[Proof of \Cref{thm: eigenvalues of K}]
	Suppose that $j= 2k$ is even. We work on $V_{k} = \{v \in L^2([0,1],\mathbb{R}^m) : v_j(0) = v_j(1) =0, \, 0<j\le k\}$. Then
	\begin{equation*}
		Q(v) = Q_{2k}(v) + R_k(v) = \int_0^1 \langle A_t v_{k}(t),v_{k}(t)\rangle dt + R_k(v)
	\end{equation*}
	
	Since the matrix $A_t$ is analytic we can diagonalize it piecewise analytically in $t$ (see \cite{kato}). Thus there exists a piecewise analytic orthogonal matrix $O_t$ such that $O_t^*A_tO_t$ is diagonal. By the second part of \Cref{prop magnitude eigenvalues}, if we make the change of coordinates $v_t\mapsto O_t v_t$ we can reduce to study the direct sum of $m$ $1-$ dimensional forms. Without loss of generality we consider forms of the type:
	\begin{equation*}
		Q_{2k}(v) = \int_0^1a_t\vert\vert v_k(t)\vert\vert^2 dt =\int_0^1a_t v_k(t)^2 dt
	\end{equation*}   
	where now $a_t$ is piecewise analytic and $v_k$ a scalar function. 
	
	For simplicity we can assume that $a_t$ does not change sign and is analytic on the whole interval. If that were not the case, we could just divide $[0,1]$ in a finite number of intervals and study $Q_{2k}$ separately on each of them. 
	
	Suppose you pick a point $t_0$ in $(0,1)$ and consider the following  subspace of codimension $mk$ in $V_k$:
	\begin{equation*}
		V_k\supset 	V^{t_0}_k = \{v \in V_k : v_j(0) = v_j(t_0) = v_j(1) =0 ,  \, 0<j\le k\}
	\end{equation*}
	
	For $t\ge t_0,$ define $v_j^{t_0}: = \int_{t_0}^tv^{t_0}_{j-1}(\tau)d\tau$ and $v_0 = v \in V_k$. It is straightforward to check that on $V_k^{t_0}$ the form $Q_{2k}$ splits as a direct sum:
	$$Q_{2k}(v) = \int_0^{t_0}\langle  A_t v_{k}(t),v_{k}(t)\rangle dt +\int_{t_0}^1\langle  A_t v^{t_0}_{k}(t),v^{t_0}_{k}(t)\rangle dt$$
	
	Now by \Cref{prop capacity} (points $i)$ and $ii)$) we can introduce as many points as we want and work separately on each segment and the asymptotic will not change (as long as the number of point is finite).
	
	Now we fix a partition $\Pi$ of $[0,1]$, $\Pi = \{t_0 =0,t_1 \dots t_{l-1},t_l =1\}$.
	Consider the subspace $V_{\Pi} = \{v \in L^2  \,\vert \, v_s(t_i)=v_s(t_{i+1})=0, 0 < s \le k, \, t_i \in \Pi\}$	 which has codimension equal to $k\vert\Pi\vert$.
	Set $a_i^- = \min_{t \in [t_i,t_{i+1}]}a_t$ and $a_i^+ = \max_{t \in [t_i,t_{i+1}]}a_t$. Finally define $v_k^{t_i}(t) = \int_{t_i}^t\dots \int_{t_i}^{\tau_1}v(\tau)d \tau \dots d \tau_{k-1} $. It follows immediately that on $V_\Pi$:
	\begin{equation*}
		\sum_i a^-_i \int_{t_i}^{t_{i+1}} v^{t_i}_k(t)^2 dt \le Q_{2k}(v) \le \sum_i a^+_i \int_{t_i}^{t_{i+1}} v^{t_i}_k(t)^2 dt
	\end{equation*}
	
	Now, we already analysed the spectrum for the problem with constant $a_t$ on $[0,1]$. The last step to understand the quantities on the right and left hand side is to see how the eigenvalues rescale when we change the length of $[0,1]$.
	
	If we look back at the proof of \Cref{lemma: boundary value problem}, it is straightforward to check that the length is relevant only when we impose the boundary conditions, we find that the eigenvalues are:
	$\lambda = \frac{a  \ell^{2k}}{(2\pi n)^{2k}}$ and again double. 
	
	In particular the estimates in \cref{eq: bound eigenvalues linear even,eq: bound eigenvalues linear odd} are still true replacing $\mu_i$ with $a_i^\pm \ell^{2k}$.
	  
	If we replace now $\ell$ by $\vert t_{i+1}-t_i \vert$ and sum the capacities according to \Cref{prop capacity} we have the following estimate on the eigenvalues on $V_\Pi$, for $n\ge 2 k \vert \Pi \vert $:
	\begin{equation*}
		\Big(\frac{\sum_i(a_i^-)^{\frac{1}{2k}}(t_{i+1}-t_i)}{\pi (n+2\vert\Pi \vert k+p(n))}\Big)^{2k} \le \lambda_n(Q_{2k}\big \vert_{V_\Pi} ) \le \Big(\frac{\sum_i(a_i^+)^{\frac{1}{2k}}(t_{i+1}-t_i)}{\pi(n-2\vert \Pi \vert k-p(n))}\Big)^{2k}
	\end{equation*}
	Moreover the min-max principle implies that, for $n\ge k\vert \Pi \vert$:
	$$\lambda_{n}\big(Q_{2k}\big \vert_{V_\Pi} \big) \le \lambda_n\big(Q_{2k}\big) \le \lambda_{n-k\vert\Pi \vert}\big(Q_{2k}\big \vert_{V_\Pi} \big)  $$
	In particular for $n\ge 3k\vert \Pi \vert$ we have:
	\begin{equation}
		\label{eq: estimate eigenvalue I proof}
		\Big(\frac{\sum_i(a_i^-)^{\frac{1}{2k}}(t_{i+1}-t_i)}{\pi (n+2\vert\Pi \vert k+p(n))}\Big)^{2k} \le \lambda_n(Q_{2k} ) \le \Big(\frac{\sum_i(a_i^+)^{\frac{1}{2k}}(t_{i+1}-t_i)}{\pi(n-3\vert \Pi \vert k-p(n))}\Big)^{2k}
	\end{equation}
	
	We address now the issue of the convergence of the Riemann sums. Set $I^\pm_a = \sum_i(a_i^\pm)^{\frac{1}{2k}}(t_{i+1}-t_i) $ and $I_a = \int_0^1a^{\frac{1}{2k}}dt$.	
	It is well know that $I^\pm_a \to I_a$ as long as $\sup_i \vert t_i-t_{i+1}\vert$ goes to zero. We need a more quantitative bound on the rate of convergence. Using results from \cite{convergenceRiemannSums} for and equispaced partition, we have that:
	\begin{equation*}
		\vert I_a-I^\pm_a\vert \le C^\pm_a \frac{1}{\vert\Pi\vert} = \frac{C(a,k,\pm)}{\mathrm{codim}(V_\Pi)}
	\end{equation*}
	Where $C(a,k,\pm)$ is a constant that depends only on the function $a$ and on $k$ and the inequality holds for $\vert\Pi\vert\ge n_0$ sufficiently large, where $n_0$ depends just on $a$ and $k$. 
	
	Consider the right hand side of \cref{eq: estimate eigenvalue I proof}, adding and subtracting $\frac{I_a}{(\pi n)^{2k}}$, we find that for $ n\ge \max\{n_0,k\vert\Pi\vert\}$:
	\begin{equation*}
		 \lambda_n(Q_{2k} ) \le \Big(\frac{I_a}{ \pi n}\Big)^{2k} + \Big(\frac{I_a^+}{\pi(n-3\vert \Pi \vert k-p(n))}\Big)^{2k}-\Big(\frac{I_a}{ \pi n}\Big)^{2k}.
	\end{equation*}
	
	A simple algebraic manipulation shows that there are constants $C_1,C_2$ and $C_3$ such that the difference on the right hand side is bounded by\begin{equation*}
		\frac{C_1 n^{2k} \vert \Pi\vert^{-1} + C_2(n^{2k}-\vert \Pi\vert^{2k}(n/\vert \Pi\vert-1)^{2k})}{ C_3(n-3k\vert \Pi\vert)^{2k} n^{2k}}
	\end{equation*}
	for $n\ge \max\{3k\vert\Pi\vert ,n_1\vert\Pi\vert,n_0\}$ where $n_1$ is a certain threshold independent of $\vert \Pi\vert$.
	
	The idea now is to choose for $n$ a partition $\Pi$ of size $\vert \Pi \vert = \lfloor n^\delta \rfloor$ to provide a good estimate of $\lambda_n(Q)$. The better result in terms of approximation is obtained for $\delta= \frac{1}{2}$. Heuristically this can be explained as follows: on one hand the first piece of the error term is of order $n^{-2k-\delta}$, comes from the convergence of the Riemann sums and gets better as $\delta\to 1$. On the other hand the second term comes from the estimate on the eigenvalues and get worse and worse as $n^\delta$ becomes comparable to $n$. 
	
	A perfectly analogous argument allows to construct an error function for the left side of \cref{eq: estimate eigenvalue I proof} which decays as $n^{-2k-1/2}$ for $n$ sufficiently large. 
	
	We have proved so far that, for one dimensional forms, $Q_{2k}$ has $2k-$capacity $\xi_+ = (\int_0^1 \sqrt[2k]{a_t}dt)^{2k}$. Now we apply point $i)$ of \Cref{prop capacity} to obtain the formula in the statement for forms on $L^2([0,1],\mathbb{R}^m)$. Finally notice that by \Cref{prop magnitude eigenvalues} the eigenvalues of $R_k(v)$ decay as $n^{-2k-1}$. If we apply point $iii)$ of \Cref{prop capacity} we find that $Q_{2k}(v)+ R_k(v)$ has the same $2k-$capacity as $Q_{2k}$ with remainder of order $1/2$.
	
	Now we consider the case $j =2k-1$. The idea is to reduce to the case of $j=4k-2$ as in the proof of \Cref{lemma: boundary value problem} and use the symmetries of $Q_{2k-1}$ to conclude.
	In the same spirit as in the beginning of the proof let us diagonalize the kernel $A_{2k-1}$. We thus reduce everything to the two dimensional case, i.e. to the quadratic forms:
	\begin{equation}
		\label{eq: odd case 2-dimensional Q}
		Q(v) = \int_0^1 \langle v_k(t),\begin{pmatrix}
			0 & -a_t \\ a_t & 0
		\end{pmatrix} v_{k-1}(t) \rangle dt \quad a_t \ge 0
	\end{equation} 
	
	It is clear that the map $v_0 \mapsto O v_0$  where $O = \begin{pmatrix}
		0 & 1 \\ 1 & 0
	\end{pmatrix}$ is an isometry of $L^{2}([0,1],\mathbb{R}^{2})$ and $Q(O v_0) = - Q(v_0)$ and so the spectrum is two sided and the asymptotic is the same for positive and negative eigenvalues.
	
	Now we reduce the problem to the even case. Let's consider \textit{the square} of $Q_{2k-1}$. By proposition \eqref{prop magnitude eigenvalues} $Q_{2k-1}$ has the same asymptotic as the form:
	\begin{equation*}
		\hat{Q}_{2k-1} = (-1)^{k+1}\int_{0}^{1} \langle A_t v_{2k-1}(t),v_0(t)\rangle dt \qquad F(v_0)(t) = (-1)^{k+1} A_tv_{2k-1}(t)
	\end{equation*}
	So we have to study the eigenvalues of the symmetric part of $F$. It is clear that:
	\begin{equation*}
			\frac{(F+F^*)^2}{4} = \frac{F^2+ F F^*+F^*F+(F^*)^2}{4}
	\end{equation*}  
	Thus we have to deal with the quadratic form:
	\begin{equation*}
		\begin{split}
			4\tilde Q(v) &= \langle [2F^{2}+F^*F + FF^*](v),v\rangle  \\ &= 2 \langle F(v),F^*(v)\rangle +  \langle F^*(v),F^*(v)\rangle + \langle F(v),F(v)\rangle
		\end{split}
	\end{equation*}
	
	The last term is the easiest to write, it is just:
	\begin{equation*}
		\langle F(v),F(v) \rangle = \int_0^1 \langle -A_t^{2} v_{2k-1}(t),v_{2k-1}(t) \rangle dt
	\end{equation*}
	which is precisely of the form of point $i)$ and gives $\frac{1}{4}$ of the desired asymptotic. 
	The operator $F^*$ acts as follows:
	\begin{equation*}
		F^*(v)=(-1)^{k+1} \int_0^t\int_0^{t_{2k-1}}\dots \int_0^{t_1} A_{t_1} v_0{(t_1)} dt_1 \dots dt_{2k-1}
	\end{equation*}
	Using integration by parts one can single out the term $A_t v_{2k-1}$. To illustrate the procedure, for $k = 1$ one gets:
	\begin{equation*}
		\begin{split}
			F^*(v) &= A_t v_1(t)-\int_0^t \dot{A}_\tau v_1(\tau) d\tau\\
			\langle F^*(v),F^*(v) \rangle &= \int_0^1 \langle -A_t^2 v_1(t),v_1(t) \rangle dt + 2\int_0^1 \langle A_t v_1(t),\int_0^t \dot{A}_\tau v_1(\tau) d\tau\rangle dt+\\
			&\qquad+\int_0^1 \langle \int_0^t \dot{A}_\tau v_1(\tau) d\tau,\int_0^t \dot{A}_\tau v_1(\tau) d\tau\rangle dt
		\end{split}
	\end{equation*}
	
	The other terms thus do not affect the asymptotic since by \Cref{prop magnitude eigenvalues} they decay at least as $O(n^{3})$. The proof goes on the same line for general $k$.
	
	The same reasoning applies to the term $\langle F(v),F^*(v)\rangle$. Summing everything one gets that the leading term is $\int_0^1 \langle -A_t^{2} v_{2k-1}(t),v_{2k-1}(t) \rangle dt$ and so this is precisely the same case as point $i)$. Recall that $A_t$ is a $2 \times 2 $ skew-symmetric matrix as defined in \cref{eq: odd case 2-dimensional Q}, thus the eigenvalues of the square coincide and are $a_t^2$. It follows that, for $n$ sufficiently large, the square of the eigenvalues of $\tilde{Q}$ satisfy:
	\begin{equation*}
		\lambda_n(\tilde Q) = \frac{\big(\int_0^1  2 \sqrt[4k-2]{a_t^2}  dt\big )^{4k-2}}{\pi^{4k-2}n^{4k-2}} + O(n^{-4k-2-\frac{1}{2}})
	\end{equation*} 

	It is immediate to see that $\frac{\big(\int_0^12 \sqrt[4k-2]{a_t^2}  dt\big)^{4k-2}}{(\pi n)^{4k-2}} = \frac{\big(\int_0^1 \sqrt[2k-1]{a_t}  dt\big)^{4k-2}}{(\pi n/2)^{4k-2}}$. This mirrors the fact that the spectrum of $Q_{2k-1}$ is double and any couple $\lambda,-\lambda$ is sent to the same eigenvalue $\lambda^2$. Thus the $(2k-1)-$capacity of $Q_{2k-1}$ is $(\int_0^1 \sqrt[2k-1]{a_t}  dt)^{2k-1}$. 
	
	Moreover, given two sequences $\{a_n\}_{n \in \mathbb{N}}$ and $\{b_n\}_{n \in \mathbb{N}}$,  $\sqrt{a_n^2+b_n^2} = a_n\sqrt{1+\frac{b_n^2}{a_n^2}} \approx a_n(1+ \frac{b_n}{a_n}+O(\frac{b_n}{a_n}))$ so the remainder is still $2k-1+\frac{1}{2}$.
	
	Arguing again by point $i)$ of \Cref{prop capacity} one gets the estimate in the statement. 
	
	The last part about the $\infty-$capacity follow just by \Cref{prop magnitude eigenvalues}. If $A_j \equiv 0$ for any $j$ then for any $\nu \in \mathbb{R}$, $\nu>0$ we have $\lambda_nn^\nu \to 0$ as $n \to \pm \infty$.
\end{proof}

		 \section{Proof of \Cref{thm: characterization of K}}
		 \label{section: proof thm 2}
		 
\begin{proof}[Proof of \Cref{thm: characterization of K}]
	The proof of the first part of the statement follows from a couple of elementary considerations. In the sequel we will use the short-hand notation $\cA$ for $Skew(K)$.
	\fact{\Cref{eq: K restricted to V is self adjoint} holds if and only if $\cA$ has finite rank}
	Suppose that $K\vert_{\mathcal{V}}$ is symmetric. Consider the orthogonal splitting of $L^2[0,1]$ as $\mathcal{V} \oplus \mathcal{V}^{\perp}$. \Cref{eq: K restricted to V is self adjoint} can be reformulated as $\cA(\mathcal{V})\subseteq \mathcal{V}^{\perp}$, thus $\imm(\cA(L^2[0,1]))\subseteq  \mathcal{V}^{\perp} + \cA( \mathcal{V}^{\perp} ) $ which is finite dimensional.
	
	Conversely, if the range of $\cA$ is finite dimensional, we can decompose $L^2[0,1]$ as $\imm(\cA)\oplus \ker(\cA)$, where the decomposition is orthogonal by skew-symmetry. Thus, on $\ker(\cA)$, $K$ is symmetric. 
	
	\fact{$\cA$ determines the kernel of $K$}
	It is well known that, if $K$ is Hilbert-Schmidt, then $K^*$ is Hilbert-Schmidt too. Since we are assuming \cref{eq: being volterra type} it is given by:
	\begin{equation*}
		K^*(v)(t) = \int_t^{1} V^*(\tau,t) v(\tau) d\tau.
	\end{equation*}
	So we can write down the integral kernel $A(t,\tau)$ of $\cA$ as follows: 
	\begin{equation*}
		A(t,\tau) = \begin{cases}
			\frac{1}{2} V(t,\tau) \text{ if } \tau <t \\
			-\frac{1}{2} V^*(\tau,t) \text{ if } t <\tau.
		\end{cases}
	\end{equation*}
	The key observation now is that the support of the kernel of $K$ is disjoint form the support of the kernel of $K^*$. Thus the kernel of $\cA$ determines the kernel of $K$ (and vice versa).
	
	Now, since we are assuming that $\cA$ has finite dimensional image, we can present its kernel as:
	\begin{equation*}
		A(t,\tau) = \frac{1}{2}  Z_t^* \cA_0 Z_\tau ,
	\end{equation*}
	where $\cA_0$ is a skew-symmetric matrix and $Z_t$ is a $ \dim(\imm(\cA)) \times k$ matrix that has as rows the elements of some orthonormal base of $\imm(\cA)$. Without loss of generality we can assume $\cA_0 = J$. In fact with an  orthogonal change of coordinates $\cA_0$ decomposes as a direct sum of rotation with an amplitude $\lambda_i$. Rescaling the coordinates by $\sqrt{\lambda_i}$ yields the desired canonical form $J$.	
	
	The first part of the statement is proved so we pass to second one. First of all notice that,
	now that we have written down any operator satisfying \cref{eq: K restricted to V is self adjoint,eq: being volterra type} in the same form as those in \cref{eq: compact part second variation}, we can apply all the results about the asymptotic of their eigenvalues. In particular, if we assume that the space $ \imm(\cA) \subset L^2([0,1],\mathbb R^k)$ is generated by piecewise analytic functions, the ordered sequence of eigenvalues satisfies:
	$$\lambda_n = \frac{\xi}{\pi n} + O(n^{-5/3}), \quad \text{ as } n \to \pm \infty.$$ 
	
	Notice that we are using a better estimates  on the reminder (for the case of the $1-$capacity) then the one given in \Cref{thm: eigenvalues of K} that was given in \cite{determinant}. We denote by $M^\dagger = \bar M^*$ the conjugate transpose. Set $2m = \dim (\imm(\cA))$, since the map $t \mapsto Z_t$ is analytic, there exists a piecewise analytic family of unitary matrices $G_t$ such that:
\begin{equation*}
	G_t^\dagger Z_t^*JZ_tG_t = \begin{bmatrix}
		&i\zeta_1(t)\\ &&\ddots\\ &&&i\zeta_l(t)\\&&&&-i\zeta_1(t)\\ &&&&&\ddots\\ &&&&&&-i\zeta_l(t)\\ &&&&&&&\underline{0}
	\end{bmatrix}
\end{equation*}
Without loss of generality we can assume that the function $\zeta_i$ are analytic on the whole interval and everywhere non negative. 
Recall that the coefficient $\xi$ appearing in the asymptotic was computed as $\xi = \int_0^1 \zeta(t)dt = \int_0^1 \sum_{i=0}^l \zeta_i(t) dt$.

Let us work on the Hilbert space $L^2([0,1],\mathbb{C}^k)$ with standard hermitian product. 
Notice that $G: L^{2}([0,1],\mathbb{C}^k )\to L^2([0,1],\mathbb{C}^k)$, $v \mapsto G_t v$ is an isometry, thus the eigenvalue of $Skew(K)=\cA$ remain the same if we consider the similar operator $G^{-1} \circ \cA  \circ G$ which acts as follows:
\begin{equation*}
	G^{-1} \circ \cA \circ G (v) = \frac{1}{2} {G}_t^\dagger Z_t^*J\int_{0}^1 Z_{\tau}G_{\tau}v(\tau) d\tau
\end{equation*}

To simplify notation let's forget about this change of coordinates and still call $Z_t$ the matrix $Z_tG_t$. Write $Z_t$ as: $$Z_t= \begin{pmatrix}
	y^*_1(t)\\ \vdots\\y_m^*(t)\\x_1^*(t)\\\vdots\\ x^*_m(t)\\
\end{pmatrix}.$$  We introduce the following notation: for a vector function $v_i$ the quantity $(v_i)_j$ stands for $j-$th component of $v_i$. 

We can now bound the function $\zeta(t)$ in terms of the components of the matrix $Z_t$:
\begin{equation*}
	\begin{split}
		2 \zeta (t) &= \sum_{j=1}^k \vert(Z_t^\dagger JZ_t)_{jj}\vert \le \sum_{i=1}^m\sum_{j=1}^{k} \vert(x_i)_j (\bar{y}_i)_j-(y_i)_j(\bar{x}_i)_j\vert(t)\\ &=\sum_{i=1}^m\sum_{j=1}^{k} 2\vert\text{Im}((x_i)_j(\bar{y}_i)_j)\vert \le \sum_{i=1}^m\sum_{j=1}^{k} 2\vert(x_i)_j\vert\vert(y_i)_j\vert = \sum_{i=1}^m 2\langle \vert x_i\vert,\vert y_i\vert \rangle(t)
	\end{split}
\end{equation*}

Where the vector $\vert v\vert$ is the vector with entries the absolute values the entries of $v$. Integrating and using H\"older inequality for the $2$ norm, we get:
\begin{equation*}
	\xi = \int_0^1 \zeta(t) dt = \sum_{i=1}^{m} \vert  \vert x_i\vert \vert_2 \, \vert\vert y_i \vert\vert_2.
\end{equation*}

The next step is to relate the quantity on the right hand side to the eigenvalues of $\cA$. The strategy now is to modify the matrix $Z_t$ in order to get an orthonormal frame of $Im(\cA)$. Keeping track of the transformations used we get a matrix representing $\cA$, then it is enough to compute the eigenvalues of the said matrix.

We can assume, without loss of generality that $\langle x_i,x_j \rangle_{L^2}  =\delta_{ij}$. This can be achieved with a symplectic change of the matrix $Z_t$. Then we modify the $y_j$ in order to make them orthogonal to the space generated by the $x_j$.  We use the following transformation:
\begin{equation*}
	\begin{pmatrix}
		Y_t \\X_t
	\end{pmatrix} \mapsto \begin{pmatrix}
	1 & M \\   0& 1\end{pmatrix} \begin{pmatrix}
	Y_t \\X_t
\end{pmatrix} = \begin{pmatrix}
Y_t+MX_t \\X_t
\end{pmatrix}
\end{equation*}
where $M$ is defined by the relation $\int_0^1 Y_tX_t^* +MX_tX_t^* dt =\int_0^1 Y_tX_t^*dt +M=0$. The last step is to make $y_j$ orthonormal. If we multiply $Y_t$ by a matrix $L$ we find the equation $L\int_0^1 Y_tY_t^* dtL^* = 1$ , so $L = (\int_0^1Y_tY_t^*dt)^{-\frac{1}{2}}$. Thus the matrix representing $\cA$ in this coordinates is one half of:
\begin{equation*}
	\cA_0 = \begin{pmatrix}
	L^{-1} &0\\
	-M^* &1
\end{pmatrix} \begin{pmatrix}
0 &-1\\ 1 &0
\end{pmatrix} \begin{pmatrix}
L^{-1} &-M \\ 0 &1
\end{pmatrix} = \begin{pmatrix}
0 &L^{-1}\\ -L^{-1}& M^*-M
\end{pmatrix} 
\end{equation*}

If we square $\cA_0$ and compute the trace we get:
$$-\frac{1}{2}\tr(\cA_0^2) = \tr(L^{-2})-\frac{1}{2} \tr ((M^*-M)^2)\ge \tr\left(\int_0^1Y_tY_t^* dt\right) = \sum_{i=1}^m \vert \vert y_i \vert \vert_2^2$$

Call $\Sigma(\cA)$ the spectrum of $\cA$, since $\cA$ is skew-symmetric it follows that: $$-\frac{1}{2}\tr(\cA_0^2) = 4 \sum_{\mu \in \Sigma(\cA), -i\mu>0 } -\mu^2 \ge0.$$

Recalling that $\vert \vert x_i \vert \vert =1$ and putting all together we find that:
\begin{equation*}
	\xi \le \sum_{i=1}^m \vert \vert y_i \vert \vert_2 \le \sqrt{m} \sqrt{\sum_{i=1}^m \vert \vert y_i \vert \vert_2^2}= 2\sqrt{m}  \sqrt{\sum_{\mu \in \Sigma(\cA), -i\mu>0 } -\mu^2}.
\end{equation*}
\end{proof}

\begin{exmpl}
	
	Consider a matrix $Z_t$ of the following form:
	\begin{equation*}
		Z_t = \begin{bmatrix}
			\xi_1(t) &\xi_3(t) \\0  &\xi_2(t)
		\end{bmatrix}
		\quad 	Z_t^*JZ_t = \begin{bmatrix}
			0 & -{\xi}_1\xi_2(t) \\ {\xi}_2\xi_1(t) & 0 
		\end{bmatrix}
	\end{equation*}
	The capacity of $K$ is given by $\zeta = \int_0^1 \vert\xi_1 \xi_2 \vert(t) dt$. We can assume that $\langle \xi_2,\xi_3 \rangle = 0$ and $\vert\vert\xi_2\vert\vert=1$. A direct computation shows that the eigenvalue of $Skew K$ are $\frac{\pm i}{2} \sqrt{(\vert\vert\xi_1\vert\vert^2+\vert\vert\xi_3\vert\vert^2)}$. This shows that the two quantities behave in a very different way. If we choose $\xi_2$ very close to $\xi_1$ and $\xi_3$ small, capacity and eigenvalue square are comparable. If we choose $\xi_3$ very big the capacity remains the same whereas the eigenvalues explode. In particular there cannot be any lower bound of $\zeta$ in terms of the eigenvalues of $K$.  
\end{exmpl}

\begin{rmrk}
	
	There is a natural class of translations that preserves the capacity. Take any path $\Phi_t$ of symplectic matrices (say $L^2$ integrable), the operators constructed with $Z_t$ and $\Phi_t Z_t$ have the same capacity (but the respective skew-symmetric part clearly do not have the same eigenvalues). 
	
	Set $K^{\Phi}(v) = \int_{0}^{t}Z_t^*J \Phi_t^{-1}\Phi_\tau Z_\tau v_\tau d \tau$ and $\Sigma^+(K^\Phi)$ the set of eigenvalues of $Skew(K^\Phi)$ satisfying $-i\sigma\ge 0$.  It seems natural to ask if:
	\begin{equation*}
		\zeta(K) = 2\inf_{\Phi_t \in Sp(n)} \sqrt{\sum_{\sigma \in \Sigma^+(K^{\Phi})}-\sigma^2}
	\end{equation*}
	
	Take for instance the example above and suppose for simplicity that $\xi_1$ and $\xi_2$ are positive and never vanishing. Using the following transformation we obtain:\begin{equation*}
		Z'_t = \begin{bmatrix}
			\sqrt{\frac{\xi_2}{\xi_1}} &\frac{-\xi_3}{\sqrt{\xi_1\xi_2}} \\ 0 & \sqrt{\frac{\xi_1}{\xi_2}}
		\end{bmatrix} \begin{bmatrix}
			\xi_1 & \xi_3 \\ 0 &\xi_2
		\end{bmatrix} = \begin{bmatrix}
			\sqrt{\xi_1\xi_2}  & 0 \\ 0 &\sqrt{\xi_1\xi_2}
		\end{bmatrix}
	\end{equation*} 
	and in this case the eigenvalue became $\frac{\pm i}{2}\langle \xi_1,\xi_2\rangle$, precisely half the capacity.
\end{rmrk}

		 \section{The second variation of an optimal control problem}
		 \label{section: control}
		We start this section collecting some basic fact about optimal control problems, first and second variation. Standard references on the topic are \cite{determinant}, \cite{bookcontrol}, \cite{bookSubriemannian}, \cite{bookJean} and \cite{symplecticMethods}.
\subsection{Symplectic geometry and optimal control problems}
Consider a smooth manifold $M$, its cotangent bundle $T^*M$ is a vector bundle on $M$ whose fibre at a point $q$ is the vector space of linear functions on $T_qM$, the tangent space of $M$ at $q$. 

Let $\pi$ be the natural projection, $\pi: T^*M \to M$ which takes a covector and gives back the base point:
\begin{equation*}
	\pi :T^*M \to M, \quad \pi(\lambda_q) = q.
\end{equation*}
Using the the projection map we define the following $1-$form, called tautological (or Liouville ) form: take an element  $X \in T_\lambda(T^*M)$, $s_\lambda(X) = \lambda(\pi_*X)$. One can check that $\sigma = ds$ is not degenerate in local coordinates. 
We obtain a symplectic manifold considering $(T^*M,\sigma)$. 

Using the symplectic form we can associate to any function on $T^*M$ a vector field.  Suppose that $H$ is a smooth function on $T^*M$, we define $\vec H$ setting:
\begin{equation*}
	\sigma(X,\vec H_\lambda ) = d_\lambda H (X), \quad \forall X \in T_\lambda(T^*M)
\end{equation*}
$H$ is called Hamiltonian function and $\vec{H}$ is an Hamiltonian vector field. 

On $T^*M$ we have a particular instance of this construction which can be used to lift arbitrary flows on the base manifold $M$ to Hamiltonian flows on $T^*M$. For any vector field $V$ on $M$ consider the following function:
\begin{equation*}
	h_V(\lambda) = \langle \lambda, V \rangle, \quad \lambda \in T^*M.
\end{equation*}
It is straight forward to check in local coordinates that $\pi_* \vec{h}_V = V$. 

The next objects we are going to introduce are Lagrangian subspaces. We say that a subspace $W$ of a symplectic vector space $(\Sigma, \sigma)$ is Lagrangian if the restriction of the symplectic form $\sigma$ is degenerate, i.e. if $\{v \in \Sigma : \sigma(v,w) =0, \, \forall \,w \in W \} = W$. An example of Lagrangian subspaces is the fibre, i.e. the kernel of $\pi_*$. More generally we can consider the following submanifolds in $T^*M$:
\begin{equation*}
	A(N)=\{\lambda \in T^*M : \lambda (X) =0, \, \forall \, X \in TN, \pi(\lambda)\in N\}
\end{equation*}
where $N\subset M$ is a submanifold. $A(N)$ is called the annihilator of $N$ and its tangent space at any point is a Lagrangian subspace.

Suppose we are given a family of complete and smooth vector fields $f_u$ which depend on some parameter $u \in U \subset \mathbb{R}^k$ and a Lagrangian, i.e. a smooth function $\varphi(u,q)$ on $U \times M$.
We use the vector fields $f_u$ to produce a family of curves on $M$. For any function $u\in L^{\infty}([0,1],U)$ we consider the following non autonomous $ODE$ system on $M$:
\begin{equation}
	\label{eq: ode admissible curve}
	\dot{q} = f_{u(t)}(q), \quad q(0) = q_0 \in M
\end{equation}
The solution are always Lipschitz curves. For fixed $q_0$, the set of functions $u \in L^{\infty}([0,1],U)$ for which said curves are defined up to time $1$ is an open set which we call $\mathcal U_{q_0}$. We can let the base point $q_0$ vary and consider $\mathcal U = \cup_{q_0\in M} \mathcal{U}_{q_0}$. It turns out that this set has a structure of a Banach manifold (see \cite{beschastnyi_morse}). We call the $L^\infty$ functions obtained this way \emph{admissible controls} and the corresponding trajectories on $M$ \emph{admissible curves}.

Denote by $\gamma_{u}$ the admissible curve obtained form  an admissible control $u$. We are interested in the following minimization problem on the space of \emph{admissible} controls:
\begin{equation}
	\label{eq: minimaztion problem}
	\min_{u \text{ admissible}} \mathcal{J}(u) = \min_{u \text{ admissible}} \int_0^1 \varphi(u(t),\gamma_{u}(t))dt
\end{equation}

We often reduce the space of admissible variations imposing additional constraints on the final and initial position of the trajectory. For example one can consider trajectories that start and end at two fixed points $q_0,q_1\in M$, or trajectory that start from a submanifold $N_0$ and reach a second submanifold $N_1$. More generally we can ask that the curves satisfy $(\gamma(0), \gamma(1))\in N\subseteq M\times M$.

We often consider the following family of functions on $T^*M$:
\begin{equation*}
	h_u: T^*M \to \mathbb R, \quad h_u(\lambda) = \langle \lambda,f_u\rangle +\nu \varphi(u,\pi(\lambda)).
\end{equation*}
We use them to lift vector fields on $M$ to vector fields on $T^*M$. They are closely relate with the function defined above and still satisfy $\pi_*(\vec{h}_u) = f_u$.

In particular, if $\tilde \gamma$ is and admissible curve, we can build a lift, i.e. a curve $\tilde\lambda$ in $T^*M$ such that $\pi(\tilde\lambda) =  \tilde \gamma$, solving $\dot \lambda = \vec h_u(\lambda)$.  
The following theorem, known as Pontryagin Maximum Principle, gives a characterization of critical points of $\mathcal{J}$, for any set of boundary conditions.

\begin{theorem*}[PMP]
	If a control $\tilde u\in L^\infty([0,1],U)$ is a critical point for the functional in \cref{eq: minimaztion problem} there exists a curve $\lambda:[0,1]\to T^*M$ and an admissible curve $q: [0,1] \to M$ such that for almost all $t\in [0,1]$
	\begin{enumerate}
		\item $\lambda(t)$ is a lift of $q(t)$: 
		$$
		q(t) = \pi (\lambda(t));
		$$
		\item $\lambda(t)$ satisfies the following Hamiltonian system:
		$$
		\frac{d \lambda}{dt} = \vec{h}_{\tilde u(t)}(\lambda);
		$$
		\item the control $\tilde u$ is determined by the maximum condition:
		$$
		h_{\tilde u(t)}(\lambda(t)) = \max_{u\in U} h_u( \lambda(t)), \quad \nu\le0;
		$$
		\item the non-triviality condition holds: $(\lambda(t),\nu)\neq (0,0)$;
		\item transversality condition holds:
		$$(-\lambda(0),\lambda(1)) \in A(N).$$
	\end{enumerate}
We call $q(t)$ an extremal curve (or trajectory) and $\lambda(t)$ an extremal.
\end{theorem*}

There are essentially two possibility for the parameter $\nu$, it can be either $0$ or, after appropriate normalization of $\lambda_t$, $-1$.
The extremals belonging to the first family are called \emph{abnormal} whereas the ones belonging to second \emph{normal}.

\subsection{The Endpoint map and its differentiation}

	We will consider now in detail the minimization problem in equation \cref{eq: minimaztion problem} with fixed endpoints.
	
	As in the previous section we denote by $\mathcal{U}_{q_0}\subset L^{\infty}([0,1],U)$ be the space of admissible controls at point $q_0$ and define the following map:
	\begin{equation*}
		E^t: \mathcal U_{q_0} \to M,\quad u\mapsto \gamma_{u}(t)	
	\end{equation*}
	It takes the control $u$ and gives the position at time $t$ of the solution of \cref{eq: ode admissible curve} starting from $q_0$. We call this map \emph{Endpoint map}. It turns out that $E^t$ is smooth, we are going now to compute its differential and Hessian.
	 The proof of these facts can be found in the book \cite{bookcontrol} or in \cite{ASZ}.
	
	For a fixed control $\tilde u $ consider the function $h_{\tilde{u}}(\lambda) = h_{\tilde{u}(t)}(\lambda)$ and define the following non autonomous flow which plays the role of parallel transport in this context:
\begin{equation}
	\label{eq: parallel transport}
	\frac{d}{dt} \tilde{\Phi}_t =  \vec{h}_{\tilde{u}}(\tilde{\Phi}_t) \qquad \tilde{\Phi}_0 = Id
\end{equation}

It has the following properties:
\begin{itemize}
	\item [\emph{i)}] It extends to the cotangent bundle the flow which solves $\dot{q} = f^t_{\tilde{u}}(q)$ on the base. In particular if $\lambda_t$ is an extremal with initial condition $\lambda_0$, $\pi(\tilde{\Phi}_t(\lambda_0)) = q_{\tilde{u}}(t)$ where $q_{\tilde{u}}$ is an extremal trajectory.
	\item [\emph{ii)}] $\tilde{\Phi}_t$ preserves the fibre over each $q \in M$. The restriction $\tilde{\Phi}_t:\,  T^*_qM \to T^*_{\tilde{\Phi}_t(q)}M $ is an affine transformation.
\end{itemize}

We suppose now that $\lambda(t)$ is an extremal and $\tilde{u}$ a critical point of the functional $\mathcal{J}$.
We use the symplectomorphism $\tilde{\Phi}_t$ to pull back the whole curve $\lambda(t)$ to the starting point $\lambda_0$. We can express all the first and second order information about the extremal using the following map and its derivatives:
\begin{equation*}
	b_u^t(\lambda) = (h_u^t-h_{\tilde{u}}^t)\circ \tilde{\Phi}_t(\lambda)
\end{equation*}

Notice that:
\begin{itemize}
	\item $b_u^t(\lambda_0)\vert_{u =\tilde{u}(t)} =0 = d_{\lambda_0}\, b_u^t\vert_{u =\tilde{u}(t)}$ by definition.
	\item $\partial_u b_u^t\vert_{u =\tilde{u}(t)} = \partial_u (h_u^t\circ \tilde{\Phi}_t)\vert_{u =\tilde{u}(t)} =0$ since $\lambda(t)$ is an extremal and $\tilde{u}$ the relative control.
\end{itemize}

Thus the first non zero derivatives are the order two ones. We define the following maps:
\begin{equation}
	\label{eq: Z_t and H_t}
	\begin{split}
		Z_t  = \partial_u \vec{b}_u^t(\lambda_0)\vert_{u=\tilde{u}(t)} : \mathbb{R}^k = T_{\tilde{u}(t)}U \to T_{\lambda_0}(T^*M) \\
		H_t = \partial_u^2 b_t(\lambda_0)\vert_{u=\tilde{u}(t)} :  \mathbb{R}^k =T_{\tilde{u}(t)}U \to  T^*_{\tilde{u}(t)}U =\mathbb{R}^k
	\end{split}
\end{equation}

We denote by $\Pi=\ker \pi_*$ the kernel of the differential of the natural projection $\pi: T^*M \to M$.
\begin{prop}[Differential of the endpoint map]
	\label{prop: differential end point map}
	Consider the endpoint map $E^t: \mathcal{U}_{q_0} \to M$. Fix a point $\tilde{u}$ and consider the symplectomorphism $\tilde{\Phi}_t$ and the map $Z_t$ defined above. The differential is the following map:
	\begin{equation*}
		d_{\tilde{u}} E (v_t) =d_{\lambda(t)} \pi \circ d_{\lambda_0}\tilde{\Phi}_t(\int_0^t Z_\tau v_\tau d\tau) \in T_{q_t}M
	\end{equation*}
\end{prop}
	In particular, if we identify $T_{\lambda_0}(T^* M)$ with $\mathbb{R}^{2m}$ and write $Z_t = \begin{pmatrix}
		Y_t \\ X_t
	\end{pmatrix}$, $\tilde{u}$ is a regular point if and only if $v_t\mapsto \int_0^t X_\tau v_\tau d\tau$ is surjective. Equivalently if the following matrix is invertible:
	\begin{equation*}
		\Gamma_t = \int_0^t X_\tau X^*_\tau d\tau \in Mat_{n\times n}(\mathbb{R}), \quad \det(\Gamma_t)\ne 0
	\end{equation*}
	
	If $d_{\tilde{u}} E^t$ is surjective then $(E^t)^{-1}(q_t)$ is smooth in a neighbourhood of $\tilde{u}$ and is tangent space is given by:
	\begin{equation*}
		\begin{split}
			T_{\tilde{u}}(E^t)^{-1}(q_t) = \{v \in L^\infty([0,1],\mathbb{R}^k) :\, \int_0^t X_\tau v_\tau d \tau =0\} \\
			= \{v \in L^\infty([0,1],\mathbb{R}^k) :\, \int_0^t Z_\tau v_\tau d \tau \in \Pi\}
		\end{split}
	\end{equation*}

When the differential of the Endpoint map is surjective a good geometric description of the situation is possible. The set of admissible control becomes smooth (at least locally) and our minimization problem can be interpreted as a constrained optimization problem. We are looking for critical points of $\mathcal{J}$ on the submanifold $\{u \in \mathcal{U} : E^t(u) = q_1\}$. 
\begin{definition}
	\label{def: strictly normal extremal}
	We say that a normal extremal $\lambda(t)$ with associated  control $\tilde{u}(t)$  is strictly normal if the differential of the endpoint map at $\tilde{u}$ is surjective. 
\end{definition}

It makes sense to go on and consider higher order optimality conditions. At critical points is well defined (i.e. independent of coordinates) the Hessian of $\mathcal J$ (or the \emph{second variation}). 
Using chronological calculus (see again \cite{bookcontrol} or \cite{ASZ}) it is possible to write the second variation of $\mathcal{J}$ on $\ker dE^t \subseteq L^{\infty}([0,1],\mathbb{R}^k)$.
\begin{prop}[Second variation]
	\label{prop: second variation}
	Suppose that $(\lambda(t),\tilde{u})$ is a strictly normal critical point of $\mathcal{J}$ with fixed initial and final point. For any $u \in L^{\infty}([0,1],\mathbb{R}^k)$ such that $\int_0^1X_tu_t dt  =0$ the second variation of $\mathcal J$ has the following expression:
	\begin{equation*}
		d^2_{\tilde{u}}\mathcal{J}(u) = -\int_0^1 \langle H_tu_t, u_t\rangle dt - \int_0^1\int_0^t \sigma ( Z_\tau u_\tau ,Z_t u_t ) d\tau dt
	\end{equation*}
	The associated bilinear form is symmetric provided that $u,v$ lie in a subspace that projects to a Lagrangian one via the map $u \mapsto \int_0^1 Z_t u_t dt$.		
	\begin{equation*}
		d^2_{\tilde{u}}\mathcal{J}(u,v) = -\int_0^1\langle H_tu_t, v_t\rangle dt - \int_0^1\int_0^t \sigma ( Z_\tau u_\tau ,Z_t v_t ) d\tau dt
	\end{equation*} 
\end{prop}

One often makes the assumption, which is customarily called \emph{strong Legendre condition}, that the matrix $H_t$ is strictly negative definite and has uniformly bounded inverse. This guarantees that the term:
	\begin{equation*}
		\int_0^1-\langle H_tu_t,v_t\rangle dt
	\end{equation*}
	is equivalent to the $L^2$ scalar product.
\begin{definition}
	\label{def: regular extremal}
	Suppose that the set $U\subset \mathbb{R}^k$ is open, we say that $(\lambda(t),\tilde{u})$ is a \emph{regular} critical point if strong Legendre condition holds along the extremal. If $H_t\le0$ but $(\lambda(t),\tilde{u})$ does not satisfy Legendre strong condition we say that $(\lambda(t),\tilde{u})$ is \emph{singular}. If $H_t \equiv 0$ we say that it is \emph{totally singular}.
\end{definition}

Even if the extremal $(\lambda(t),\tilde{u})$ is abnormal or not strictly normal it is possible to produce a second variation for the optimal control problem. To do so one considers the extended control system:
\begin{equation*}
	\hat{f}_{(v,u)}(q) = \begin{pmatrix}
		\varphi(u,q)+v \\
		f_u(q)
	\end{pmatrix} \in \mathbb{R}\times T_qM
\end{equation*}
and the corresponding endpoint map $\hat E^t: (0,+\infty) \times \mathcal U_{q_0} \to \mathbb{R}\times M$. To differentiate it we use the same construction explained above and employ the following Hamiltonians on $\mathbb{R}^*\times T^*M$:
\begin{equation*}
	\hat{h}_{(v,u)} (\nu,\lambda)= \langle \lambda,f_u\rangle +\nu(\varphi(u,q)+v)
\end{equation*} 
One has just to identify which are the right controls to consider, PMP implies that $\dot \nu =0$, $\nu \le0$ and $v =0$. In the end one obtains formally the same expression as in \Cref{prop: second variation} involving the derivatives of the functions $\hat h_{(v,u)}$ and recover the same expression as in \Cref{prop: second variation} for strictly normal extremals (see \cite[Chapter 20]{bookcontrol} or \cite{symplecticMethods}).

\subsection{Reformulation of the main results}
In this section we reformulate \Cref{thm: characterization of K} as a characterization of the compact part of the second variation of an optimal control problem at a strictly normal regular extremal (see \cref{def: strictly normal extremal,def: regular extremal}).

\begin{thm}
	\label{theorem 1}
	Suppose $\cV\subset L^2([0,1],\mathbb{R}^k)$ is a finite codimension subspace and $K$ and operator satisfying \cref{eq: K restricted to V is self adjoint,eq: being volterra type}. Then $(K,\cV)$ can be realized as the second variation of an optimal control problem at a strictly normal regular extremal. 	
	To any such couple we can associate a triple $((\Sigma,\sigma),\Pi,Z)$ consisting of:
	\begin{itemize}
		\item  a finite dimensional symplectic space $(\Sigma,\sigma)$;
		\item   a Lagrangian subspace $\Pi\subset \Sigma$;
		\item  a linear map $Z: L^2([0,1],\mathbb{R}^k) \to \Sigma$ such that $\imm(Z)$ is transversal to the subspace $\Pi$.
	\end{itemize}
	This triple is unique up to the action of $\mathrm{stab}_\Pi (\Sigma,\sigma)$, the group of symplectic transformations that fix $\Pi$.  Any other triple is given by $((\Sigma,\sigma),\Pi,\Phi \circ Z)$ for $\Phi\in \mathrm{stab}_\Pi (\Sigma,\sigma) $.
	
	Vice versa any triple $((\Sigma,\sigma),\Pi,Z)$ as above determines a couple $(K,\cV)$. We can define the skew-symmetric part $\cA$ of $K$ as:
	\begin{equation*}
		 \langle \cA u, v \rangle = \sigma (Zu,Zv), \, \forall u,v \in L^2([0,1],\mathbb{R}^k),
	\end{equation*}
	$\cA$ determines the whole operator $K$ and its domain is recovered as $ \cV = Z^{-1}(\Pi)$.
	\begin{proof}
		The proof is essentially a reformulation of \Cref{thm: characterization of K}. Given the operator we construct the symplectic space $(\Sigma,\sigma)$ taking as vector space the image of the skew-symmetric part $\imm(\cA)$ and as symplectic form $\langle \cA \cdot, \cdot \rangle$.
		
		The transversality condition correspond to the fact that the differential of the endpoint map is surjective. 
		
		The only thing left to show is uniqueness of the triple. Without loss of generality we can assume that the symplectic subspace $(\Sigma, \sigma) = (\mathbb{R}^{2n},\sigma)$ is the standard one and that the Lagrangian subspace $\Pi$ is the vertical subspace. In this coordinates $$Z(v) = \int_0^1Z_t v_t dt  = \int_0^1 \begin{pmatrix}
			Y_t \\ X_t
		\end{pmatrix} v_t dt.$$
	
		Define the following map: 
		$$F: L^2([0,1],\mathrm{Mat}_{n\times k}(\mathbb{R})) \to  L^2([0,1]^2,\mathrm{Mat}_{k\times k}(\mathbb{R})), \quad Y_t \mapsto Z_t^*JZ_{\tau} = X_t^*Y_\tau-Y_t^*X_\tau.$$		
		It is linear  if $X_t$ is fixed. To determine uniqueness we have to study an affine equation thus is sufficient to study the kernel of $F$. Suppose for simplicity that $X_t$ and $Y_t$ are continuous in $t$. We have to solve the equation:
		\begin{equation*}
			F(Y_{t}) = Z_t^*JZ_{\tau} =\sigma (Z_t ,Z_{\tau}) = 0.
		\end{equation*}
		
		Consider the following subspace of $\mathbb{R}^{2n}$ 
		\begin{equation*}
			V^{[0,1]} = \Big \{\sum_{i=1}^l  Z_{t_i}\nu_i :  \, \nu_i \in \mathbb{R}^k, t_i \in[0,1], l \in \mathbb{N}\Big \} \subset \mathbb{R}^{2n}
		\end{equation*}
		
		It follows that $F(Y_t) =0$ if and only if the subspace $V^{[0,1]}$ is isotropic. Since we are in finite dimension, we can consider a finite number of instants $t_i$ to which we can restrict to generate the whole $V^{[0,1]}$. Call $I$ the set of this instants.
		Without loss of generality we can assume that $\{\sum_{i \in I}X_{t_i}\nu_i, \nu_i \in \mathbb{R}^{k}, t_i \in I\}=\mathbb{R}^{n}$. 
		
		This is so since the image of $Z$ is transversal to $\Pi$ and thus $\Gamma = \int_0^1 X_tX_t^* dt $ is non degenerate. In fact, if the subspace $\{\sum_{i=1}^l  X_{t_i}\nu_i \vert  \, \nu_i \in \mathbb{R}^k, l \in \mathbb{N}\}$ were a proper subspace of $\mathbb{R}^n$, there would be a vector $\mu$ such that $\langle\mu, X_t \nu\rangle  =0$, $\forall t \in [0,1]$ and $\forall \nu \in \mathbb{R}^n$. Thus an element of the kernel of $\Gamma$. A contradiction.
		
		Now we evaluate the equation $F(Y_t) =0 \iff  Y_t^* X_\tau=X_t^*Y_\tau $ at the instants $t = t_i$ that guarantee controllability.  One can read off the following identities:
		\begin{equation*}
			Y_t^* v_j = X_t^* c_j 
		\end{equation*}
		where the $v_j'$s are a base of $\mathbb{R}^n$ and $c_j$ free parameters.  Taking transpose we get that $Y_t = G X_t$.
		
		It is straightforward to check that, if $Y_t = G X_t$, $G$ must be symmetric, in fact:
		\begin{equation*}
			Z_tJZ_{\tau} =Y_t^*X_{\tau}-X_t^* Y_{\tau}=X_t^*(G^*-G)X_{\tau} = 0 \iff  G = G^*
		\end{equation*} 
		And so uniqueness is proved when $X_t$ and $Y_t$ are continuous. 
		
		The case in which $X_t$ and $Y_t$ are just $L^2$ (matrix-)functions can be dealt with similarly. One has just to replace \emph{evaluations} with integrals of the form $\int_{t-\epsilon}^{t+\epsilon}Z_\tau  \nu d\tau$ and $\int_{t-\epsilon}^{t+\epsilon} X_\tau \nu d\tau$ and interpret every equality $t$ almost everywhere. 
		
		The only thing left to show is how to construct a control system  with given $(K,\mathcal{V})$ as second variation. By the equivalence stated above it is enough to show that we can realize any given map $Z : L^2([0,1],\mathbb{R}^k)\to \Sigma$ with a proper control system. We can assume without loss of generality that $(\Sigma,\sigma)$ is just $\mathbb{R}^{2m}$ with the standard symplectic form and $\Pi$ is the vertical subspace. With this choices the map $Z$ is given by :
		\begin{equation*}
			v \mapsto \int_0^1 Z_t v_t dt =   \int_0^1 \begin{pmatrix}
				Y_t v_t  \\X_t v_t 
			\end{pmatrix} dt
		\end{equation*}		
		
		The operator $K$ is then given by $K(v) = \int_0^t Z_t^*J Z_\tau v_\tau d \tau$ and $\mathcal{V} = \{v \vert \int_0^1 X_t v_t dt  =0\}$. Consider the following linear quadratic system on $\mathbb{R}^m$: 
		\begin{equation*}
			f_u(q) = B_t u \quad \varphi_t(x) =\frac{1}{2}\vert u \vert^2+ \langle \Omega_t u,x \rangle,
		\end{equation*}
		where $B_t$ and $\Omega_t$ are matrices of size $m \times k$, the Hamiltonian in PMP reads:
		\begin{equation*}
			h_u(\lambda,x) = \langle \lambda, B_t u \rangle -\frac{1}{2}\vert u \vert^2 -\langle \Omega_t u,x \rangle
		\end{equation*}
		Take as extremal control $u_t \equiv 0$, it easy to check that the re-parametrization flow $\tilde \Phi_t$ defined in \cref{eq: parallel transport} is just the identity and the matrix $Z_t$ for this problem is the following:
		\begin{equation*}
			Z_t = \begin{pmatrix}
				\Omega_t \\B_t
			\end{pmatrix}
		\end{equation*}
		So it is enough to take $\Omega_t = Y_t$ and $B_t = X_t $.
	\end{proof}
\end{thm}


We can reformulate also the second part of \Cref{thm: characterization of K} relating the capacity of $K$ and the eigenvalues of $\cA$.
We make the following assumptions:
\begin{enumerate}
	\item  the map $t\mapsto Z_t $ is piecewise analytic in $t$;
	\item the maximum condition in the statement of PMP defines a $C^2$ function $	\hat{H}_t(\lambda) = \max_{u\in \mathbb R^k} h^t_u( \lambda)$ in a neighbourhood of the strictly normal regular extremal we are considering.
\end{enumerate}

Under the above assumptions the following proposition clarifies the link between the matrices $Z_t$ and $ H_t$ and the function $\hat{H}_t$. A proof can be found either in \cite[Proposition 21.3]{bookcontrol} or \cite{ASZ}.
\begin{prop}
		Suppose that $(\lambda(t),\tilde{u})$ is an extremal and the function $\hat H_t$ is $C^2$, using the flow defined in \cref{eq: parallel transport}  define $\mathcal{H}_t(\lambda) = (\hat H_t-h_{\tilde u(t)})\circ \tilde \Phi_t(\lambda) $. It holds that:
		\begin{equation*}
		 \text{Hess}_{\lambda_0}(\mathcal H_t) = JZ_tH_t^{-1}Z_t^*J 
		\end{equation*}
	\end{prop}

Define $R_t =  \max_{v \in \mathbb{R}^k,\vert\vert v\vert\vert=1} \vert\vert Z_t v\vert\vert$ and let $\{\pm i\zeta_j(t)\}_{j=1}^l$ be the eigenvalues of $iZ_t^*JZ_t$ as defined in \Cref{section: proof thm 2}. We have the following proposition.
\begin{prop}
	\label{prop maximized hamiltonian capacity}
	The capacity $\xi$ of $K$ satisfies:
	\begin{equation*}
		\xi \le \frac{\sqrt{k} \, \vert\vert R_t\vert\vert_2}{2}  \sqrt{\int_{0}^{1} \tr( \text{Hess}_{\lambda_0}(\mathcal H_t)) dt}
	\end{equation*} and in particular, if we order the functions $\zeta_j(t)$ decreasingly, they satisfy \begin{equation*}
		0 \le \zeta_j(t) \le R_t \sqrt{\lambda_{2j}(t)}, \quad j \in \{1,\dots l\}
	\end{equation*} 
	where $\lambda_{j}(t)$ are the eigenvalues of $Hess_{\lambda_0}(\mathcal H_t)$ in decreasing order. 
	\begin{proof}
	We give a sketch of the proof. Without loss of generality we can assume $H_t =- Id$, otherwise, we can perform the change of coordinate on $L^2([0,1],\mathbb{R}^k)$ $v \mapsto (-H_t)^{-\frac{1}{2}}v$ and redefine $Z_t$ accordingly.
	
	In this notation $Hess_{\lambda_0}(\mathcal{H}_t)$ corresponds to the matrix $JZ_tZ_t^*J$. If we square $A_t = Z_t^*JZ_t$ we obtain:
	\begin{equation*}
		A_t^*A_t = -Z_t^*JZ_tZ_t^*JZ_t = - Z_t^*\big (J Z_tZ_t^* J\big)Z_t = -Z_t^*Hess_{\lambda_0}(\mathcal{H}_t)Z_t
	\end{equation*}
	Observe that $\zeta_j(t)$ is an eigenvalue of $A_t$ if and only if $-\zeta_j^2(t)$ is a eigenvalue of $A^*_tA_t$. The equation above relates the \emph{restriction} of $Hess_{\lambda_0}(\mathcal{H}_t)$ to the image of the maps $Z_t: \mathbb{R}^k\to \mathbb{R}^{2n}$ with the square of the functions $\zeta_j(t)$ defining the capacity.
	
	The idea is to use Cauchy interlacing inequality for the eigenvalues of $Hess_{\lambda_0}(\mathcal{H}_t)$ and its restriction to a codimension $2n-k$ subspace. If $\{\lambda_j(t)\}_{j=1}^{2n}$ are the eigenvalues of the Hessian, taken in decreasing order, and $\{\mu_j(t)\}_{j=1}^{2n-k}$  the eigenvalues of its restriction we have:
	$$\lambda_{j+2n-k}(t) \le \mu_j(t) \le \lambda_j(t) $$ 
	In our case $Z_t$ are not orthogonal projectors but we can adjust the estimates considering how much the matrices $Z_t$ dilate the space, and thus we have to take in account the function $R_t$ defined just before the statement. Denote by $\mu_j(t)$ the $j-$th  eigenvalue of $-A_t^2,$ putting all together we have:
	$$0 \le \mu_j(t) \le R_t^2 \, \lambda_{2j}(t)\quad j \in \{1,\dots k\}$$
	Where we shifted the index by one since $\mu_{2k-1}(t) = \mu_{2k}(t)$ for all $k\le l$. Taking square roots and integrating we have:
	$$\int_0^1\zeta_j(t)dt\le \int_0^1R_t \sqrt{\lambda_{2j}(t)}dt $$ 
	Summing up over $j$ we find that:
	$$\xi = \int_0^1\sum_j \zeta_j(t) dt \le \frac{1}{2}\int_0^1 \sum_j R_t\sqrt{\lambda_{2j}(t)}dt \le \frac{\sqrt{k}\vert \vert R_t \vert \vert_2}{2} \sqrt{\int_0^1 \tr (Hess_{\lambda_0} (\mathcal{H}_t))}$$
	\end{proof}
\end{prop}

We turn now to \Cref{thm: eigenvalues of K}, we can interpret it as a quantitative version of various necessary optimality conditions that one can formulate for certain classes of singular extremals (see \cite[Chapter 20]{bookcontrol} or \cite[Chapter 12]{bookSubriemannian}). Moreover, leaving optimality conditions aside, \Cref{thm: eigenvalues of K} gives the asymptotic distribution of the eigenvalues of the second variation for totally singular extremals (see \cref{def: regular extremal}).

As mentioned in the previous section we can produce a second variation also in the non strictly normal case which is at least formally very similar to the normal case.
However, a common occurrence is that the matrix $H_t$ completely degenerates and is constantly equal to the zero matrix. This is the case for affine control systems and abnormal extremal in Sub-Riemannian geometry, i.e. systems of the form:
\begin{equation*}
	f_u = \sum_{i=1}^l f_i u_i +f_0, \quad f_i \text{ smooth vector fields}
\end{equation*}

In this case Legendre condition $H_t\le0$ (see the previous section) does not give much information. One, then, looks for \emph{higher} order optimality conditions. This is usually done exactly as in \Cref{lemma taylor Q}: the first optimality conditions one finds are \emph{Goh condition} and  \emph{ generalized Legendre condition} which prevent the second variation from being \emph{strongly indefinite}. 

In the notation of \Cref{lemma taylor Q} Goh conditions is written as $Q_1 \equiv 0$ i.e. $Z_t^*JZ_t \equiv 0$.
It can be reformulated in geometric terms as follows, if $\lambda_t$ is the extremal then
$$\lambda_t [\partial_u f_u(q(t))v_1,\partial_u f_u(q(t))v_2] =0,\, \forall \, v_1,v_2 \in \mathbb{R}^k$$
From \Cref{thm: eigenvalues of K} it is clear that if $Q_1\not \equiv 0$, the second variation has infinite negative index and that eigenvalues distribute evenly between the negative and positive parts of the spectrum. Then one asks that the second term $Q_2$ is non positive definite (recall the different sign convention in \Cref{prop: second variation}), otherwise the negative part of the spectrum of $-Q_2$ becomes infinite. In our notation this  condition reads $$(Z_t^{(1)})^*JZ_t \le 0 \iff \sigma (Z_t^{(1)} v, Z_t v) \le 0, \, \forall \, v  \in \mathbb{R}^k.$$ Again it can be translated in a differential condition along the extremal, however this time it will in general involve more than just commutators if the system is not control affine.

If $Q_2\equiv 0$, one can take more derivatives and find new conditions. In particular, using the notation of \Cref{lemma taylor Q}, one has always to ask that the first non zero term in the expansion is of even order and that the matrix of its coefficients is non positive  in order to have finite negative index.

		\section*{Acknowledgements}
		The author wishes to thank Prof. A. Agrachev for the stimulating discussions on the topic and the referee for the helpful suggestions which greatly improved the exposition.
		
		\bibliographystyle{plain}
		\bibliography{ref}		
\end{document}